\documentclass{article}
\usepackage{amsfonts}
\usepackage{graphicx}
\usepackage{amsmath}
\usepackage{amssymb}
\newtheorem{theorem}{Theorem}

\newtheorem{corollary}[theorem]{Corollary}

\newtheorem{definition}[theorem]{Definition}
\newtheorem{example}[theorem]{Example}

\newtheorem{lemma}[theorem]{Lemma}

\newtheorem{observation}[theorem]{Observation}

\newtheorem{proposition}[theorem]{Proposition}
\newtheorem{remark}[theorem]{Remark}

\newenvironment{proof}[1][Proof]{\textbf{#1.} }{\ \rule{0.5em}{0.5em}}

\begin{document}

\title{Classification of $\left(n-5\right)$-filiform Lie algebras\thanks{%
Research partially supported by the D.G.I.C.Y.T project PB98-0758}}
\author{Jos\'{e} Mar\'{\i}a Ancochea Berm\'{u}dez\thanks{%
corresponding author: Tel +00 34 913944566: fax: 00 34 91 3944564, e-mail: Jose\_Ancochea@mat.ucm.es} \and Otto Rutwig
Campoamor Stursberg \\
Departamento de Geometr\'{\i}a y Topolog\'{\i}a\\
Fac. CC. Matem\'{a}ticas Univ. Complutense\\
28040 Madrid ( Spain )}
\date{}
\maketitle

\begin{abstract}
In this paper we consider the problem of classifying the $(n-5)$-filiform
Lie algebras. This is the first index for which infinite parametrized
families appear, as can be seen in dimension $7.$ Moreover we obtain large
families of characteristic nilpotent Lie algebras with nilpotence index $5$
and show that at least for dimension $10$ there is a characteristic
nilpotent Lie algebra with nilpotence index $4$ which is the algebra of
derivations of a nilpotent Lie algebra.

\textit{Keywords: }$p$-filiform, characteristically nilpotent; Lie algebras
\end{abstract}

\section{\textbf{Generalities}}

Nilpotent Lie algebras have played an important role in mathematics over the
last thirty years: either in the classification theory of Lie algebras,
where they play a central role as a consequence of the L\'{e}vi theorem and
the works of Malcev, or in the geometrical and analytical applications such
as the nilmanifolds, which allow to construct concrete compact differential
manifolds, or Pfaffian systems.\newline
The first important research about nilpotent Lie algebras is due to K.Umlauf
in the last 19$^{th}$ century. In the 40's and 50's Morozov and Dixmier
begun with the systematical study of this class of algebras. Morozov gave a
classification of six dimensional nilpotent Lie algebras in $1958$ [16]. The
existence of an infinity of complex nilpotent Lie algebras from dimension
seven up showed the complexity of the classification problem. A complete
classification of $7$-dimensional nilpotent Lie algebras was obtained by the
first author and M.Goze [3].\newline
We pointed out that for dimensions greater or equal than eight only partial
classifications are known. The most of them correspond to the filiform Lie
algebras, i.e. algebras with maximal nilpotence index. They are classified
up to dimension $11$ [10]. It seems natural to determine an invariant which
measures the nilpotence of Lie algebras. The first author and M.Goze
introduced in [2] an invariant that allowed not only the
classification in dimension seven, but the study of nilpotent Lie algebras
with lower nilpotence indexes.\newline
Let $\frak{g}_{n}=\left( \Bbb{C}^{n},\mu _{n}\right) $ be a nilpotent Lie
algebra. For each $X\in \Bbb{C}^{n}$ we denote $c\left( X\right) $ the
ordered sequence of a similitude invariant of the nilpotent operator $ad_{%
\frak{g}_{n}}\left( X\right) ,$ i.e. the ordered sequence of dimensions of
the Jordan blocks for this operator. We consider the lexicographical order
in the set of these sequences.

\begin{definition}
The characteristic sequence of $\frak{g}_{n}$ is an isomorphism invariant $%
c\left( \frak{g}_{n}\right) $ defined by 
\[
c\left( \frak{g}_{n}\right) =\max_{X\in \frak{g}_{n}-C^{1}\frak{g}%
_{n}}\left\{ c\left( X\right) \right\} 
\]
where $C^{1}\frak{g}_{n}$ is the derived algebra. A nonzero vector $X\in 
\frak{g}_{n}-C^{1}\frak{g}_{n}$ satisfying $c\left( X\right) =c\left( \frak{g%
}_{n}\right) $ is called characteristic vector.
\end{definition}

\begin{definition}
A nilpotent Lie algebra $\frak{g}_{n}$ is called $p$-filiform if its
characteristic sequence is $c\left( \frak{g}_{n}\right) =\left(
n-p,1,..^{\left( p\right) }..,..,1\right) .$
\end{definition}

\begin{remark}
It follows immediately from the definition of $p$-filiformness that the $%
\left( n-1\right) $-filiform Lie algebras are the abelian algebras $\frak{a}$%
. It is easily shown that the $\left( n-2\right) $-filiform Lie algebras are
the direct sum of an Heisenberg algebra $\frak{H}_{2p+1}$ and an abelian
algebra. A classification of the $\left( n-3\right)$- and $\left( n-4\right)$-filiform Lie algebras can be found in [7], respectively [6]. 
The former is also the last where the
number of isomorphism classes is finite, as we shall see. We are primarly interested in nonsplit $\left( n-5\right) $-filiform Lie
algebras, for the general $\left( n-5\right) $-filiform algebras are
obtained by direct sums of nonsplit algebras and abelian algebras. Thus
classifying the nonsplit we have classified all of them.
\end{remark}

\section{\textbf{The classification theorem}}

\begin{theorem}
$\left( \text{ Classification theorem}\right) $\newline
Each $n$-dimensional nonsplit $\left( n-5\right) $-filiform Lie algebra $%
\frak{g}_{n}$ is isomorphic to one of the laws $\mu _{n}^{i},\;i\in \left\{
1,..,103\right\} $ listed below.
\end{theorem}

Before we give the list in even and odd dimension we have to introduce some
notation. This will be applicable to both odd and even dimensions. Let $%
\frak{g}_{n}$ be a $n$-dimensional nilpotent complex Lie algebra. Then we
identify the Lie algebra with its law $\left( \Bbb{C}^{n},\mu _{n}\right) ,$
where $\mu _{n}\in \frak{T}_{\left( 2,1\right) }^{a}$ is an alternated
tensor of type $\left( 2,1\right) $ satisfying the Jacobi equation. We
denote the derived subalgebra as $C^{1}\frak{g}_{n}.$ The list is structured
as follows: at first we indicate indexes for which common brackets are
listed. The bullet item completes the corresponding law. This kind of
presentation has two advantages: on one hand side it is easier to read the
concrete algebra laws, on the other it indicates in a certain manner that
the laws are closely related ( as it follows in the proof). Finally, as
usual the nonwritten brackets are zero or obtained by antisymmetry.

\subsection{Even dimension}
The first subdivision is referred to the dimension of the derived algebra. 
\bigskip
1. $\dim C^{1}\frak{g}_{n}=6$\newline
There are four laws with these conditions:\newline
For the indexes $i\in \{1,2,3\}$ we have the common brackets\newline
$\mu _{2m}^{i}\left( X_{1},X_{j}\right) =X_{j+1},\;j\in \{2,3,4,5\},\;m\geq
4 $\newline
$\mu _{2m}^{i}\left( X_{5},X_{2}\right) =\mu _{2m}^{i}\left(
X_{3},X_{4}\right) =Y_{1};\;$ $\mu _{2m}^{i}\left( Y_{2t-1},Y_{2t}\right)
=X_{6},\;2\leq t\leq m-3$ if $m>4$

\begin{itemize}
\item  $\mu _{2m}^{1}\left( X_{3},X_{2}\right) =Y_{2};\;\mu _{2m}^{1}\left(
Y_{2},X_{3}\right) =X_{6};\;\mu _{2m}^{1}\left( Y_{2},X_{2}\right) =X_{5};$

\item  $\mu _{2m}^{2}\left( X_{3},X_{2}\right) =Y_{2};\;$

\item  $\mu _{2m}^{3}\left( X_{4},X_{2}\right) =X_{6};\;\mu _{2m}^{3}\left(
X_{3},X_{2}\right) =Y_{2}+X_{5};$\newline
\end{itemize}

\bigskip
For $i=4$ we obtain the law\newline
$\mu _{2m}^{4}\left( X_{1},X_{j}\right) =X_{j+1},\;j\in \{2,3,4,5\},\;m\geq
5 $\newline
$\mu _{2m}^{4}\left( X_{5},X_{2}\right) =\mu _{2m}^{4}\left(
X_{3},X_{4}\right) =Y_{1};\;\mu _{2m}^{4}\left( X_{3},X_{2}\right) =Y_{2};$%
\newline
$\mu _{2m}^{4}\left( Y_{3},X_{3}\right) =X_{6};\;\mu _{2m}^{4}\left(
Y_{3},X_{2}\right) =X_{5};$\newline
$\mu _{2m}^{4}\left( Y_{2t-1},Y_{2t}\right) =X_{6},\;2\leq t\leq m-3$ $.$

\bigskip
2. $\dim C^{1}\frak{g}_{n}=5$

\bigskip
For $i\in \left\{ 5,6,7^{\alpha },8,...,19\right\} $ we have\newline
$\mu _{2m}^{i}\left( X_{1},X_{j}\right) =X_{j+1},\;j\in \{2,3,4,5\},\;m\geq
4 $\newline
$\mu _{2m}^{i}\left( Y_{2t-1},Y_{2t}\right) =X_{6},\;2\leq t\leq m-3$ if $%
m>4.$

\begin{itemize}
\item  $\mu _{2m}^{5}\left( X_{5},X_{2}\right) =\mu _{2m}^{5}\left(
X_{3},X_{4}\right) =Y_{1};\;\mu _{2m}^{5}\left( Y_{2},X_{3}\right) =X_{6};$%
\newline
$\mu _{2m}^{5}\left( Y_{2},X_{2}\right) =X_{5}$.

\item  $\mu _{2m}^{6}\left( X_{3},X_{2}\right) =Y_{1};\;\mu _{2m}^{6}\left(
Y_{1},X_{3}\right) =X_{6};\;\mu _{2m}^{6}\left( Y_{1},X_{2}\right)
=X_{5}+X_{6};$\newline
$\mu _{2m}^{6}\left( Y_{2},X_{2}\right) =X_{6}.$

\item  $\mu _{2m}^{7,\alpha }\left( X_{4},X_{2}\right) =\alpha X_{6};\;\mu
_{2m}^{7,\alpha }\left( X_{3},X_{2}\right) =Y_{1}+\alpha X_{5},\;\alpha \neq
0;$\newline
$\mu _{2m}^{7,\alpha }\left( Y_{1},X_{3}\right) =X_{6};\;\mu _{2m}^{7,\alpha
}\left( Y_{1},X_{2}\right) =X_{5}+X_{6};\;\mu _{2m}^{7,\alpha }\left(
Y_{2},X_{2}\right) =X_{6}.$

\item  $\mu _{2m}^{8}\left( X_{3},X_{2}\right) =Y_{1};\;\mu _{2m}^{8}\left(
Y_{1},X_{3}\right) =X_{6};\;\mu _{2m}^{8}\left( Y_{1},X_{2}\right) =X_{5};$%
\newline
$\mu _{2m}^{8}\left( Y_{2},X_{2}\right) =X_{6}.$

\item  $\mu _{2m}^{9}\left( X_{4},X_{2}\right) =X_{6};\;\mu _{2m}^{9}\left(
X_{3},X_{2}\right) =Y_{1}+X_{5};\;\mu _{2m}^{9}\left( Y_{1},X_{3}\right)
=X_{6};$\newline
$\mu _{2m}^{9}\left( Y_{1},X_{2}\right) =X_{5};\;\mu _{2m}^{9}\left(
Y_{2},X_{2}\right) =X_{6}.$

\item  $\mu _{2m}^{10}\left( X_{3},X_{2}\right) =Y_{1};\;\mu
_{2m}^{10}\left( Y_{1},X_{2}\right) =X_{6};\;\mu _{2m}^{10}\left(
Y_{2},X_{3}\right) =X_{6};$\newline
$\mu _{2m}^{10}\left( Y_{2},X_{2}\right) =X_{5}.$

\item  $\mu _{2m}^{11}\left( X_{5},X_{2}\right) =\mu _{2m}^{11}\left(
X_{3},X_{4}\right) =X_{6};\;\mu _{2m}^{11}\left( X_{3},X_{2}\right) =Y_{1};$%
\newline
$\mu _{2m}^{11}\left( Y_{1},X_{2}\right) =X_{6};\;\mu _{2m}^{11}\left(
Y_{2},X_{3}\right) =X_{6};\;\mu _{2m}^{11}\left( Y_{2},X_{2}\right) =X_{5}.$

\item  $\mu _{2m}^{12}\left( X_{3},X_{2}\right) =Y_{1};\;\mu
_{2m}^{12}\left( Y_{2},X_{3}\right) =X_{6};\;\mu _{2m}^{12}\left(
Y_{2},X_{2}\right) =X_{5};\;$

\item  $\mu _{2m}^{13}\left( X_{5},X_{2}\right) =\mu _{2m}^{13}\left(
X_{3},X_{4}\right) =X_{6};\;\mu _{2m}^{13}\left( X_{3},X_{2}\right) =Y_{1};$%
\newline
$\mu _{2m}^{13}\left( Y_{2},X_{3}\right) =X_{6};\;\mu _{2m}^{13}\left(
Y_{2},X_{2}\right) =X_{5}.\;$

\item  $\mu _{2m}^{14}\left( X_{5},X_{2}\right) =\mu _{2m}^{14}\left(
X_{3},X_{4}\right) =\;\mu _{2m}^{14}\left( X_{4},X_{2}\right) =X_{6};$%
\newline
$\mu _{2m}^{14}\left( X_{3},X_{2}\right) =Y_{1}+X_{5};\;\mu _{2m}^{14}\left(
Y_{2},X_{3}\right) =X_{6};\;\mu _{2m}^{14}\left( Y_{2},X_{2}\right) =X_{5}.$

\item  $\mu _{2m}^{15}\left( X_{4},X_{2}\right) =X_{6};\;\mu
_{2m}^{15}\left( X_{3},X_{2}\right) =Y_{1}+X_{5};\;\mu _{2m}^{15}\left(
Y_{2},X_{3}\right) =X_{6};$\newline
$\mu _{2m}^{15}\left( Y_{2},X_{2}\right) =X_{5}.$

\item  $\mu _{2m}^{16}\left( X_{3},X_{2}\right) =Y_{1};\;\mu
_{2m}^{16}\left( Y_{2},X_{2}\right) =X_{6}.$

\item  $\mu _{2m}^{17}\left( X_{5},X_{2}\right) =\mu _{2m}^{17}\left(
X_{3},X_{4}\right) =X_{6};\;\mu _{2m}^{17}\left( X_{3},X_{2}\right) =Y_{1};$%
\newline
$\mu _{2m}^{17}\left( Y_{2},X_{2}\right) =X_{6}.$

\item  $\mu _{2m}^{18}\left( X_{5},X_{2}\right) =\mu _{2m}^{18}\left(
X_{3},X_{4}\right) =\mu _{2m}^{18}\left( X_{4},X_{2}\right) =X_{6};$\newline
$\mu _{2m}^{18}\left( X_{3},X_{2}\right) =Y_{1}+X_{5};\;\mu _{2m}^{18}\left(
Y_{2},X_{2}\right) =X_{6}.$

\item  $\mu _{2m}^{19}\left( X_{4},X_{2}\right) =X_{6};\;\mu
_{2m}^{19}\left( X_{3},X_{2}\right) =Y_{1}+X_{5};\;\mu _{2m}^{19}\left(
Y_{2},X_{2}\right) =X_{6}.$
\end{itemize}

\bigskip
For $i\in \{20,21,22,23\}$\newline
$\mu _{2m}^{i}\left( X_{1},X_{j}\right) =X_{j+1},\;j\in \{2,3,4,5\},\;m\geq
5 $\newline
$\mu _{2m}^{i}\left( Y_{2},X_{3}\right) =X_{6};\;\mu _{2m}^{i}\left(
Y_{2},X_{2}\right) =X_{5};\;\mu _{2m}^{i}\left( Y_{3},X_{2}\right) =X_{6};$%
\newline
$\mu _{2m}^{i}\left( Y_{2},Y_{4}\right) =X_{6};\;\;$\newline
$\mu _{2m}^{i}\left( Y_{2t+1},Y_{2t+2}\right) =X_{6};\;2\leq t\leq m-4$ if $%
m>5.$

\begin{itemize}
\item  $\mu _{2m}^{20}\left( X_{3},X_{2}\right) =Y_{1};\;$

\item  $\mu _{2m}^{21}\left( X_{5},X_{2}\right) =\mu _{2m}^{21}\left(
X_{3},X_{4}\right) =X_{6};\;\mu _{2m}^{21}\left( X_{3},X_{2}\right) =Y_{1}.$

\item  $\mu _{2m}^{22}\left( X_{5},X_{2}\right) =\mu _{2m}^{22}\left(
X_{3},X_{4}\right) =\mu _{2m}^{22}\left( X_{4},X_{2}\right) =X_{6};$\newline
$\mu _{2m}^{22}\left( X_{3},X_{2}\right) =Y_{1}+X_{5}.$

\item  $\mu _{2m}^{23}\left( X_{4},X_{2}\right) =X_{6};\;\mu
_{2m}^{23}\left( X_{3},X_{2}\right) =Y_{1}+X_{5}.$
\end{itemize}

\bigskip
3. $\dim C^{1}\frak{g}_{n}=4$

\bigskip
For the indexes $i\in \{24,..,29\}$ the derived algebra $C^{1}\frak{g}_{n}$
is not abelian.\newline
For $i\in \{24,25\}$ we have\newline
$\mu _{2m}^{i}\left( X_{1},X_{j}\right) =X_{j+1},\;j\in \{2,3,4,5\},\;m\geq
4 $\newline
$\mu _{2m}^{i}\left( X_{5},X_{2}\right) =\mu _{2m}^{i}\left(
X_{3},X_{4}\right) =X_{6};\;$\newline
$\mu _{2m}^{i}\left( Y_{1},X_{3}\right) =X_{6};\;\mu _{2m}^{i}\left(
Y_{1},X_{2}\right) =X_{5};\;\mu _{2m}^{i}\left( Y_{2},X_{2}\right) =X_{6};$%
\newline
$\mu _{2m}^{i}\left( Y_{2t-1},Y_{2t}\right) =X_{6},\;2\leq t\leq m-3$ if $%
m>4 $

\begin{itemize}
\item  $\mu _{2m}^{25}\left( X_{4},X_{2}\right) =X_{6};\;\mu
_{2m}^{25}\left( X_{3},X_{2}\right) =X_{5}.$
\end{itemize}

\bigskip
For $i\in \{26,27\}$\bigskip
$\mu _{2m}^{i}\left( X_{1},X_{j}\right) =X_{j+1},\;j\in \{2,3,4,5\},\;m\geq
4 $\newline
$\mu _{2m}^{i}\left( X_{5},X_{2}\right) =\mu _{2m}^{i}\left(
X_{3},X_{4}\right) =X_{6};\;$\newline
$\mu _{2m}^{i}\left( Y_{2t-1},Y_{2t}\right) =X_{6},\;1\leq t\leq m-3$

\begin{itemize}
\item  $\mu _{2m}^{26}\left( Y_{1},X_{3}\right) =X_{6};\;\mu
_{2m}^{26}\left( Y_{1},X_{2}\right) =X_{5}.$

\item  $\mu _{2m}^{27}\left( X_{4},X_{2}\right) =X_{6};\;\mu
_{2m}^{27}\left( X_{3},X_{2}\right) =X_{5};\;\mu _{2m}^{27}\left(
Y_{1},X_{3}\right) =X_{6};$\newline
$\mu _{2m}^{27}\left( Y_{1},X_{2}\right) =X_{5}.$
\end{itemize}

\bigskip
For $i\in \{28,29\}$\newline
$\mu _{2m}^{i}\left( X_{1},X_{j}\right) =X_{j+1},\;j\in \{2,3,4,5\},\;m\geq
3 $\newline
$\mu _{2m}^{i}\left( X_{5},X_{2}\right) =\mu _{2m}^{i}\left(
X_{3},X_{4}\right) =X_{6};\;$\newline
$\mu _{2m}^{i}\left( Y_{2t-1},Y_{2t}\right) =X_{6},\;1\leq t\leq m-3$ if $%
m>3 $

\begin{itemize}
\item  $\mu _{2m}^{29}\left( X_{4},X_{2}\right) =X_{6};\;\mu
_{2m}^{29}\left( X_{3},X_{2}\right) =X_{5}.$
\end{itemize}

\bigskip
For the indexes $i\in \{30,..,53\}$ the derived algebra $C^{1}\frak{g}_{n}$
is abelian.\newline
For $i\in \left\{ 30,31,32,33,34,35,36\right\} $\newline
$\mu _{2m}^{i}\left( X_{1},X_{j}\right) =X_{j+1},\;j\in \{2,3,4,5\},\;m\geq
5 $\newline
$\mu _{2m}^{i}\left( Y_{1},X_{j}\right) =X_{j+2},\;j\in \{2,3,4\};\;$\newline
$\mu _{2m}^{i}\left( Y_{2},X_{j}\right) =X_{j+3},\;j\in \{2,3\};\;$\newline
$\mu _{2m}^{i}\left( Y_{2t+1},Y_{2t+2}\right) =X_{6},\;2\leq t\leq m-3$ if $%
m>5$

\begin{itemize}
\item  $\mu _{2m}^{30}\left( Y_{3},X_{2}\right) =X_{6};\;\mu
_{2m}^{30}\left( Y_{1},Y_{4}\right) =X_{6}.$

\item  $\mu _{2m}^{31}\left( Y_{3},X_{2}\right) =X_{6};\;\mu
_{2m}^{31}\left( Y_{1},Y_{4}\right) =X_{6};\;\mu _{2m}^{31}\left(
Y_{2},Y_{4}\right) =X_{6}.$

\item  $\mu _{2m}^{32}\left( Y_{3},X_{2}\right) =X_{6};\;\mu
_{2m}^{32}\left( Y_{1},Y_{4}\right) =X_{6};\;\mu _{2m}^{32}\left(
Y_{2},Y_{3}\right) =X_{6}.$

\item  $\mu _{2m}^{33}\left( Y_{3},X_{2}\right) =X_{6};\;\mu
_{2m}^{33}\left( Y_{1},Y_{4}\right) =X_{6};\;\mu _{2m}^{33}\left(
Y_{2},Y_{3}\right) =X_{6};$\newline
$\mu _{2m}^{33}\left( Y_{2},Y_{4}\right) =X_{6}.$

\item  $\mu _{2m}^{34}\left( Y_{3},X_{2}\right) =X_{6};\;\mu
_{2m}^{34}\left( Y_{2},Y_{4}\right) =X_{6}.$

\item  $\mu _{2m}^{35}\left( Y_{3},X_{2}\right) =X_{6};\;\mu
_{2m}^{35}\left( Y_{1},Y_{3}\right) =X_{6};\;\mu _{2m}^{35}\left(
Y_{2},Y_{4}\right) =X_{6}.$

\item  $\mu _{2m}^{36}\left( Y_{1},Y_{3}\right) =X_{6};\;\mu
_{2m}^{36}\left( Y_{2},Y_{4}\right) =X_{6}.$
\end{itemize}

\bigskip
For $i\in \{37,38,39,40,41\}$\newline
$\mu _{2m}^{i}\left( X_{1},X_{j}\right) =X_{j+1},\;j\in \{2,3,4,5\},\;m\geq
4 $\newline
$\mu _{2m}^{i}\left( Y_{1},X_{j}\right) =X_{j+2},\;j\in \{2,3,4\};\;$\newline
$\mu _{2m}^{i}\left( Y_{2t-1},Y_{2t}\right) =X_{6},\;1\leq t\leq m-3.$

\begin{itemize}
\item  $\mu _{2m}^{37}\left( Y_{2},X_{3}\right) =X_{6};\;\mu
_{2m}^{37}\left( Y_{2},X_{2}\right) =X_{5}.$

\item  $\mu _{2m}^{38}\left( Y_{2},X_{2}\right) =X_{6}.$

\item  $\mu _{2m}^{39}\left( X_{3},X_{2}\right) =X_{6};\;\mu
_{2m}^{39}\left( Y_{2},X_{2}\right) =X_{6}.$

\item  $\mu _{2m}^{41}\left( X_{3},X_{2}\right) =X_{6}.$
\end{itemize}

\bigskip
For $i\in \{42,43,44\}$\newline
$\mu _{2m}^{i}\left( X_{1},X_{j}\right) =X_{j+1},\;j\in \{2,3,4,5\},\;m\geq
4 $\newline
$\mu _{2m}^{i}\left( Y_{1},X_{j}\right) =X_{j+2},\;j\in \{2,3,4\}\;$\newline
$\mu _{2m}^{i}\left( Y_{2t-1},Y_{2t}\right) =X_{6},\;2\leq t\leq m-3$ if $%
m>4 $

\begin{itemize}
\item  $\mu _{2m}^{42}\left( Y_{2},X_{3}\right) =X_{6};\;\mu
_{2m}^{42}\left( Y_{2},X_{2}\right) =X_{5}.$

\item  $\mu _{2m}^{43}\left( Y_{2},X_{2}\right) =X_{6}.$

\item  $\mu _{2m}^{44}\left( X_{3},X_{2}\right) =X_{6};\;\mu
_{2m}^{44}\left( Y_{2},X_{2}\right) =X_{6}.$
\end{itemize}

\bigskip
For $i\in \left\{ 45,46\right\} $\newline
$\mu _{2m}^{i}\left( X_{1},X_{j}\right) =X_{j+2},\;j\in \{2,3,4,5\},\;m\geq
4 $\newline
$\mu _{2m}^{i}\left( Y_{1},X_{3}\right) =X_{6};\;\mu _{2m}^{i}\left(
Y_{1},X_{2}\right) =X_{5};\;\mu _{2m}^{i}\left( Y_{2},X_{2}\right) =X_{6};$%
\newline
$\mu _{2m}^{i}\left( Y_{2t-1},Y_{2t}\right) =X_{6},\;2\leq t\leq m-3$ if $%
m>4 $

\begin{itemize}
\item  $\mu _{2m}^{46}\left( X_{4},X_{2}\right) =X_{6};\;\mu
_{2m}^{46}\left( X_{3},X_{2}\right) =X_{5}.$
\end{itemize}

\bigskip
For $i\in \{47,48,49,50\}$\newline
$\mu _{2m}^{i}\left( X_{1},X_{j}\right) =X_{j+1},\;j\in \{2,3,4,5\},\;m\geq
4 $\newline
$\mu _{2m}^{i}\left( Y_{1},X_{3}\right) =X_{6};\;\mu _{2m}^{i}\left(
Y_{1},X_{2}\right) =X_{5};$\newline
$\mu _{2m}^{i}\left( Y_{2t-1},Y_{2t}\right) =X_{6},\;1\leq t\leq m-3.$

\begin{itemize}
\item  $\mu _{2m}^{47}\left( Y_{2},X_{2}\right) =X_{6}.$

\item  $\mu _{2m}^{48}\left( X_{4},X_{2}\right) =X_{6};\;\mu
_{2m}^{48}\left( X_{3},X_{2}\right) =X_{5};\;\mu _{2m}^{48}\left(
Y_{2},X_{2}\right) =X_{6}.$

\item  $\mu _{2m}^{50}\left( X_{4},X_{2}\right) =X_{6};\;\mu
_{2m}^{50}\left( X_{3},X_{2}\right) =X_{5}.$
\end{itemize}

\bigskip
For $i\in \{51,52,53\}$\newline
$\mu _{2m}^{i}\left( X_{1},X_{j}\right) =X_{j+1},\;j\in \{2,3,4,5\},\;m\geq
3 $\newline
$\mu _{2m}^{i}\left( Y_{2t-1},Y_{2t}\right) =X_{6},\;1\leq t\leq m-3$ if $%
m>3 $

\begin{itemize}
\item  $\mu _{2m}^{52}\left( X_{3},X_{2}\right) =X_{6}.$

\item  $\mu _{2m}^{53}\left( X_{4},X_{2}\right) =X_{6};\;\mu
_{2m}^{53}\left( X_{3},X_{2}\right) =X_{5}.$
\end{itemize}

\subsection{Odd dimension}

\bigskip
1. $\dim C^{1}\frak{g}_{n}=6$\newline
There is only one law:\newline
$\mu _{2m+1}^{54}\left( X_{1},X_{j}\right) =X_{j+1},\;j\in \left\{
2,3,4,5\right\} ,\;m\geq 4$\newline
$\mu _{2m+1}^{54}\left( X_{5},X_{2}\right) =\mu _{2m+1}^{54}\left(
X_{3},X_{4}\right) =Y_{1};\;\mu _{2m+1}^{54}\left( X_{3},X_{2}\right)
=Y_{2}; $\newline
$\mu _{2m+1}^{54}\left( Y_{3},X_{3}\right) =X_{6};\;\mu _{2m+1}^{54}\left(
Y_{3},X_{2}\right) =X_{5};\;\mu _{2m+1}^{54}\left( Y_{2t},Y_{2t+1}\right)
=X_{6},$\newline
$2\leq t\leq m-3$ if $m>4.$

\bigskip
2. $\dim C^{1}\frak{g}_{n}=5$\newline
For $i\in \{55,56,57,58,59,60,61\}$\bigskip
$\mu _{2m+1}^{i}\left( X_{1},X_{j}\right) =X_{j+1},\;j\in
\{2,3,4,5\},\;m\geq 4$\newline
$\mu _{2m+1}^{i}\left( Y_{2},X_{3}\right) =X_{6};\;\mu _{2m+1}^{i}\left(
Y_{2},X_{2}\right) =X_{5};\;$\newline
$\mu _{2m+1}^{i}\left( Y_{2t},Y_{2t+1}\right) =X_{6},\;1\leq t\leq m-3$

\begin{itemize}
\item  $\mu _{2m+1}^{55}\left( X_{5},X_{2}\right) =\mu _{2m+1}^{55}\left(
X_{3},X_{4}\right) =Y_{1}.$

\item  $\mu _{2m+1}^{56}\left( X_{3},X_{2}\right) =Y_{1};\;\mu
_{2m+1}^{56}\left( Y_{1},X_{2}\right) =X_{6}.$

\item  $\mu _{2m+1}^{57}\left( X_{5},X_{2}\right) =\mu _{2m+1}^{57}\left(
X_{3},X_{4}\right) =X_{6};\;\mu _{2m+1}^{57}\left( X_{3},X_{2}\right) =Y_{1};
$\newline
$\mu _{2m+1}^{57}\left( Y_{1},X_{2}\right) =X_{6}.$

\item  $\mu _{2m+1}^{58}\left( X_{3},X_{2}\right) =Y_{1}.$

\item  $\mu _{2m+1}^{59}\left( X_{5},X_{2}\right) =\mu _{2m+1}^{59}\left(
X_{3},X_{4}\right) =X_{6};\;\mu _{2m+1}^{59}\left( X_{3},X_{2}\right) =Y_{1}.
$

\item  $\mu _{2m+1}^{60}\left( X_{5},X_{2}\right) =\mu _{2m+1}^{60}\left(
X_{3},X_{4}\right) =\mu _{2m+1}^{60}\left( X_{4},X_{2}\right) =X_{6};$%
\newline
$\mu _{2m+1}^{60}\left( X_{3},X_{2}\right) =Y_{1}+X_{5}.$
\end{itemize}

\bigskip
For $i\in \{62,..,74\}$\newline
$\mu _{2m+1}^{i}\left( X_{1},X_{j}\right) =X_{j+1},\;j\in
\{2,3,4,5\},\;m\geq 3$\newline
$\mu _{2m+1}^{i}\left( Y_{2t},Y_{2t+1}\right) =X_{6},\;1\leq t\leq m-3$ if $%
m>3$

\begin{itemize}
\item  $\mu _{2m+1}^{62}\left( X_{5},X_{2}\right) =\mu _{2m+1}^{62}\left(
X_{3},X_{4}\right) =Y_{1};\;\mu _{2m+1}^{62}\left( X_{4},X_{2}\right) =X_{6};
$\newline
$\mu _{2m+1}^{62}\left( X_{3},X_{2}\right) =X_{5}.$

\item  $\mu _{2m+1}^{63}\left( X_{5},X_{2}\right) =\mu _{2m+1}^{63}\left(
X_{3},X_{4}\right) =Y_{1};\;\mu _{2m+1}^{63}\left( X_{3},X_{2}\right) =X_{6}.
$

\item  $\mu _{2m+1}^{64}\left( X_{5},X_{2}\right) =\mu _{2m+1}^{64}\left(
X_{3},X_{4}\right) =Y_{1}.$

\item  $\mu _{2m+1}^{65}\left( X_{3},X_{2}\right) =Y_{1};\;\mu
_{2m+1}^{65}\left( Y_{1},X_{3}\right) =X_{6};\;\mu _{2m+1}^{65}\left(
Y_{1},X_{2}\right) =X_{5}+X_{6}.$

\item  $\mu _{2m+1}^{66,\alpha }\left( X_{4},X_{2}\right) =\alpha
X_{6};\;\mu _{2m+1}^{66,\alpha }\left( X_{3},X_{2}\right) =Y_{1}+\alpha
X_{5},\;\alpha \neq 0;$\newline
$\mu _{2m+1}^{66,\alpha }\left( Y_{1},X_{3}\right) =X_{6};$\ $\mu
_{2m+1}^{66,\alpha }\left( Y_{1},X_{2}\right) =X_{5}+X_{6}.$

\item  $\mu _{2m+1}^{67}\left( X_{3},X_{2}\right) =Y_{1};\;\mu
_{2m+1}^{67}\left( Y_{1},X_{3}\right) =X_{6};\;\mu _{2m+1}^{67}\left(
Y_{1},X_{2}\right) =X_{5}.$

\item  $\mu _{2m+1}^{68}\left( X_{4},X_{2}\right) =X_{6};\;\mu
_{2m+1}^{68}\left( X_{3},X_{2}\right) =Y_{1}+X_{5};$\newline
$\,\mu _{2m+1}^{68}\left( Y_{1},X_{3}\right) =X_{6};\;\mu _{2m+1}^{68}\left(
Y_{1},X_{2}\right) =X_{5}.$

\item  $\mu _{2m+1}^{69}\left( X_{3},X_{2}\right) =Y_{1};\;\mu
_{2m+1}^{68}\left( Y_{1},X_{2}\right) =X_{6}.$

\item  $\mu _{2m+1}^{70}\left( X_{5},X_{2}\right) =\mu _{2m+1}^{70}\left(
X_{3},X_{4}\right) =X_{6};\;\mu _{2m+1}^{70}\left( X_{3},X_{2}\right) =Y_{1};
$\newline
$\mu _{2m+1}^{70}\left( Y_{1},X_{2}\right) =X_{6}.$

\item  $\mu _{2m+1}^{71}\left( X_{3},X_{2}\right) =Y_{1}.$

\item  $\mu _{2m+1}^{72}\left( X_{5},X_{2}\right) =\mu _{2m+1}^{72}\left(
X_{3},X_{4}\right) =X_{6};\;\mu _{2m+1}^{72}\left( X_{3},X_{2}\right) =Y_{1}.
$

\item  $\mu _{2m+1}^{73}\left( X_{5},X_{2}\right) =\mu _{2m+1}^{73}\left(
X_{3},X_{4}\right) =\mu _{2m+1}^{73}\left( X_{4},X_{2}\right) =X_{6};$%
\newline
$\mu _{2m+1}^{73}\left( X_{3},X_{2}\right) =Y_{1}+X_{5}.$

\item  $\mu _{2m+1}^{74}\left( X_{4},X_{2}\right) =X_{6};\;\mu
_{2m+1}^{74}\left( X_{3},X_{2}\right) =Y_{1}+X_{5}.$
\end{itemize}

\bigskip
For $i\in \{75,76,77,78\}$\newline
$\mu _{2m+1}^{i}\left( X_{1},X_{j}\right) =X_{j+1},\;j\in
\{2,3,4,5\},\;m\geq 4$\newline
$\mu _{2m+1}^{i}\left( Y_{2},X_{3}\right) =X_{6};\;\mu _{2m+1}^{i}\left(
Y_{2},X_{2}\right) =X_{5};\;\mu _{2m+1}^{i}\left( Y_{3},X_{2}\right) =X_{6};$%
\newline
$\mu _{2m+1}^{i}\left( Y_{2t},Y_{2t+1}\right) =X_{6},\;2\leq t\leq m-3$ if $%
m>4.$

\begin{itemize}
\item  $\mu _{2m+1}^{75}\left( X_{3},X_{2}\right) =Y_{1}.$

\item  $\mu _{2m+1}^{76}\left( X_{5},X_{2}\right) =\mu _{2m+1}^{76}\left(
X_{3},X_{4}\right) =X_{6};\;\mu _{2m+1}^{76}\left( X_{3},X_{2}\right) =Y_{1}.
$

\item  $\mu _{2m+1}^{77}\left( X_{5},X_{2}\right) =\mu _{2m+1}^{77}\left(
X_{3},X_{4}\right) =\mu _{2m+1}^{77}\left( X_{4},X_{2}\right) =X_{6};$%
\newline
$\mu _{2m+1}^{77}\left( X_{3},X_{2}\right) =Y_{1}+X_{5}.$

\item  $\mu _{2m+1}^{78}\left( X_{4},X_{2}\right) =X_{6};\;\mu
_{2m+1}^{78}\left( X_{3},X_{2}\right) =Y_{1}+X_{5}.$
\end{itemize}

\bigskip
3. $\dim C^{1}\frak{g}_{n}=4$\newline
For the indexes $i\in \{79,..,85\}$ the derived algebra $C^{1}\frak{g}_{n}$
is not abelian.

\bigskip
For $i\in \{79,80\}$\newline
$\mu _{2m+1}^{i}\left( X_{1},X_{j}\right) =X_{j+1},\;j\in
\{2,3,4,5\},\;m\geq 4$\newline
$\mu _{2m+1}^{i}\left( X_{5},X_{2}\right) =\mu _{2m+1}^{i}\left(
X_{3},X_{4}\right) =X_{6};\;$\newline
$\mu _{2m+1}^{i}\left( Y_{1},X_{j}\right) =X_{j+3},\;j\in \{2,3\};$\newline
$\mu _{2m+1}^{i}\left( Y_{2},X_{2}\right) =\mu _{2m+1}^{i}\left(
Y_{1},Y_{3}\right) =X_{6};$\newline
$\mu _{2m+1}^{i}\left( Y_{2t},Y_{2t+1}\right) =X_{6},\;2\leq t\leq m-3$ if $%
m>4.$

\begin{itemize}
\item  $\mu _{2m+1}^{80}\left( X_{4},X_{2}\right) =X_{6};\;\mu
_{2m+1}^{80}\left( X_{3},X_{2}\right) =X_{5}.$
\end{itemize}

\bigskip
For $i\in \{81,82,83,84,85\}$\newline
$\mu _{2m+1}^{i}\left( X_{1},X_{j}\right) =X_{j+1},\;j\in
\{2,3,4,5\},\;m\geq 3$\newline
$\mu _{2m+1}^{i}\left( X_{5},X_{2}\right) =\mu _{2m+1}^{i}\left(
X_{3},X_{4}\right) =X_{6};\;$\newline
$\mu _{2m+1}^{i}\left( Y_{2t},Y_{2t+1}\right) =X_{6},\;1\leq t\leq m-3$ if $%
m>3.$

\begin{itemize}
\item  $\mu _{2m+1}^{81}\left( Y_{1},X_{3}\right) =X_{6};\;\mu
_{2m+1}^{81}\left( Y_{1},X_{2}\right) =X_{5}+X_{6}.$

\item  $\mu _{2m+1}^{82}\left( Y_{1},X_{3}\right) =X_{6};\;\mu
_{2m+1}^{82}\left( Y_{1},X_{2}\right) =X_{5}.$

\item  $\mu _{2m+1}^{83}\left( X_{4},X_{2}\right) =X_{6};\;\mu
_{2m+1}^{83}\left( X_{3},X_{2}\right) =X_{5};\;\mu _{2m+1}^{83}\left(
Y_{1},X_{3}\right) =X_{6};$\newline
$\mu _{2m+1}^{83}\left( Y_{1},X_{2}\right) =X_{5}.$

\item  $\mu _{2m+1}^{84}\left( Y_{1},X_{2}\right) =X_{6}.$

\item  $\mu _{2m+1}^{85}\left( X_{4},X_{2}\right) =X_{6};\;\mu
_{2m+1}^{85}\left( X_{3},X_{2}\right) =X_{5};\;\mu _{2m+1}^{85}\left(
Y_{1},X_{2}\right) =X_{6}$.
\end{itemize}

\bigskip
For the indexes $i\in \{86,..,103\}$ the derived algebra $C^{1}\frak{g}_{n}$
is abelian.

\bigskip
For $i=86$ we have \newline
$\mu _{2m+1}^{86}\left( X_{1},X_{j}\right) =X_{j+1},\;j\in
\{2,3,4,5\},\;m\geq 5$\newline
$\mu _{2m+1}^{86}\left( Y_{1},X_{j}\right) =X_{j+2},\;j\in \{2,3,4\}\;$%
\newline
$\mu _{2m+1}^{86}\left( Y_{2},X_{j}\right) =X_{j+3},\;j\in \{2,3\};$\newline
$\mu _{2m+1}^{86}\left( Y_{3},X_{2}\right) =\mu _{2m+1}^{86}\left(
Y_{1},Y_{4}\right) =X_{6};$\newline
$\mu _{2m+1}^{86}\left( Y_{2},Y_{5}\right) =X_{6};$\ $\mu _{2m+1}^{86}\left(
Y_{2t},Y_{2t+1}\right) =X_{6},\;3\leq t\leq m-3$ if $m>5.$

\bigskip
For $i\in \{87,..,92\}$\newline
$\mu _{2m+1}^{i}\left( X_{1},X_{j}\right) =X_{j+1},\;j\in
\{2,3,4,5\},\;m\geq 4$\newline
$\mu _{2m+1}^{i}\left( Y_{1},X_{j}\right) =X_{j+2},\;j\in \{2,3,4\};$\newline
$\mu _{2m+1}^{i}\left( Y_{2},X_{j}\right) =X_{j+3},\;j\in \{2,3\}\;\;$%
\newline
$\mu _{2m+1}^{i}\left( Y_{2t},Y_{2t+1}\right) =X_{6},\;2\leq t\leq m-3$ if $%
m>4.$

\begin{itemize}
\item  $\mu _{2m+1}^{87}\left( Y_{3},X_{2}\right) =\mu _{2m}^{87}\left(
Y_{2},Y_{3}\right) =X_{6}.$

\item  $\mu _{2m+1}^{88}\left( Y_{3},X_{2}\right) =X_{6}.$

\item  $\mu _{2m+1}^{89}\left( Y_{3},X_{2}\right) =\mu _{2m+1}^{89}\left(
Y_{1},Y_{2}\right) =X_{6}.$

\item  $\mu _{2m+1}^{90}\left( Y_{3},X_{2}\right) =\mu _{2m+1}^{90}\left(
Y_{1},Y_{3}\right) =X_{6}.$

\item  $\mu _{2m+1}^{91}\left( Y_{1},Y_{3}\right) =X_{6}.$

\item  $\mu _{2m+1}^{92}\left( Y_{2},Y_{3}\right) =X_{6}.$
\end{itemize}

\bigskip
For $i\in \{93,94,95\}$\newline
$\mu _{2m+1}^{i}\left( X_{1},X_{j}\right) =X_{j+1},\;j\in
\{2,3,4,5\},\;m\geq 3$\newline
$\mu _{2m+1}^{i}\left( Y_{1},X_{j}\right) =X_{j+2},\;j\in \{3,4\}$\newline
$\mu _{2m+1}^{i}\left( Y_{2t},Y_{2t+1}\right) =X_{6},\;1\leq t\leq m-3$ if $%
m>3.$

\begin{itemize}
\item  $\mu _{2m+1}^{93}\left( Y_{1},X_{2}\right) =X_{4}+X_{6}.$

\item  $\mu _{2m+1}^{94}\left( Y_{1},X_{2}\right) =X_{4}.$

\item  $\mu _{2m+1}^{95}\left( Y_{1},X_{2}\right) =X_{4};\;\mu
_{2m+1}^{95}\left( X_{3},X_{2}\right) =X_{6}.$
\end{itemize}

\bigskip
For $i\in \{96,97\}$\newline
$\mu _{2m+1}^{i}\left( X_{1},X_{j}\right) =X_{j+1},\;j\in
\{2,3,4,5\},\;m\geq 4$\newline
$\mu _{2m+1}^{i}\left( Y_{1},X_{j}\right) =X_{j+3},\;j\in \{2,3\};$\newline
$\mu _{2m+1}^{i}\left( Y_{2},X_{2}\right) =\mu _{2m+1}^{i}\left(
Y_{1},Y_{3}\right) =X_{6};$\newline
$\mu _{2m+1}^{i}\left( Y_{2t},Y_{2t+1}\right) =X_{6},\;2\leq t\leq m-3$ if $%
m>4$

\begin{itemize}
\item  $\mu _{2m+1}^{97}\left( X_{4},X_{2}\right) =X_{6};\;\mu
_{2m+1}^{97}\left( X_{3},X_{2}\right) =X_{5}.$
\end{itemize}

\bigskip
For $i\in \{98,..,103\}$\newline
$\mu _{2m+1}^{i}\left( X_{1},X_{j}\right) =X_{j+1},\;j\in
\{2,3,4,5\},\;m\geq 3$\newline
$\mu _{2m+1}^{i}\left( Y_{2t},Y_{2t+1}\right) =X_{6},\;1\leq t\leq m-3$ if $%
m>3$

\begin{itemize}
\item  $\mu _{2m+1}^{98}\left( X_{4},X_{2}\right) =X_{6};\;\mu
_{2m+1}^{98}\left( X_{3},X_{2}\right) =X_{5};\;\mu _{2m+1}^{98}\left(
Y_{1},X_{3}\right) =X_{6};$\newline
$\mu _{2m+1}^{98}\left( Y_{1},X_{2}\right) =X_{5}.$

\item  $\mu _{2m+1}^{99}\left( Y_{1},X_{3}\right) =X_{6};\;\mu
_{2m+1}^{99}\left( Y_{1},X_{2}\right) =X_{5}.$

\item  $\mu _{2m+1}^{100}\left( Y_{1},X_{3}\right) =X_{6};\;\mu
_{2m+1}^{100}\left( Y_{1},X_{2}\right) =X_{5}+X_{6}.$

\item  $\mu _{2m+1}^{101}\left( Y_{1},X_{2}\right) =X_{6}.$

\item  $\mu _{2m+1}^{102}\left( X_{4},X_{2}\right) =X_{6};\;\mu
_{2m+1}^{102}\left( X_{3},X_{2}\right) =X_{5};\;\mu _{2m+1}^{102}\left(
Y_{1},X_{2}\right) =X_{6}.$

\item  $\mu _{2m+1}^{103}\left( X_{3},X_{2}\right) =\mu _{2m+1}^{103}\left(
Y_{1},X_{2}\right) =X_{6}.$
\end{itemize}

\section{\textbf{Distinction}}

We now prove that the Lie algebras listed above are pairwise non isomorphic.
For this purpose we establish the following notation: We denote by $\frak{l}%
_{2m}^{j}$ and $\frak{l}_{2m+1}^{j}$ the $j^{th}$ class of $\left(
n-4\right) $-filiform nonsplit Lie algebras in even, respectively odd
dimension. See [6] for this classification. From the nonsplitness of the
algebras listed above it follows immediately that the quotient through the
center will be a $\left( n-4-k\right) $-filiform Lie algebra in dimension $%
\left( n-k\right) $ for the values $k=1,2,3.$ We will call quotient type to
the nonsplit part of $\frac{\frak{g}_{n}}{Z\left( \frak{g}_{n}\right) }$ \
for $\frak{g}_{n}$ from the list, where $Z\left( \frak{g}_{n}\right) $ is
the center of the algebra, and write $\left[ \frac{\frak{g}_{n}}{Z\left( 
\frak{g}_{n}\right) }\right] $.\newline
For certain algebras having the same quotient type it will be useful to
consider ideal classes with specific assumptions. A $7$-dimensional
nilpotent complex Lie algebras is noted $\frak{n}_{k}^{7}$, where $k$
denotes the isomorphism class. For the classification of these algebras see
[3].

\begin{lemma}
Let $\frak{g}_{n}^{i}\;$be a Lie algebra from the list above.

\begin{enumerate}
\item  if $n=2m$ and $i\in \left\{ 6,7,..,23\right\} $ $\left( \text{or }%
n=2m+1\text{ and }i\in \left\{ 55,56,..,78\right\} \right) $there exists a
unique class $[\frak{n}]$ of $7$-dimensional ideals with characteristic
sequence $\left( 5,1,1\right) .$

\item  if $n=2m$ and $i\in \{51,52,53\}\;\left( \text{or }n=2m+1\text{ and }%
i\in \{101,102,103\}\right) $ there exists a unique class $[\frak{N}]$ of $6$%
-dimensional filiform ideals.%
\endproof%
\end{enumerate}
\end{lemma}

It will be always necessary to distinguish isomorphism classes by the
dimension of its algebra of derivations $Der\left( \frak{g}_{n}\right) $ or
the weights of diagonalizable derivations of the algebra $\frak{g}_{n}.$ The
characteristic polynomial is noted $p_{c}\left( \lambda \right) $ with $%
\lambda $ the variable. A derivation $f\in Der\left( \frak{g}_{n}\right) $
is written: 
\[
f\left( X_{i}\right)
=\sum_{j=1}^{6}f_{i}^{j}X_{j}+\sum_{k=1}^{n-6}g_{i}^{k}Y_{k},\;\;f\left(
Y_{i}\right) =\sum_{j=1}^{6}h_{i}^{j}X_{j}+\sum_{k=1}^{n-6}l_{i}^{j}Y_{j} 
\]
for $f_{i}^{j},g_{i}^{j},h_{i}^{j},l_{i}^{j}\in \Bbb{C}.\;\;$

\begin{observation}
For those Lie algebras having the same dimension for the algebra of
derivations the distinction will follow from the multiplicities of the
weights for diagonalizable derivations. This information is also comprised
in tables $8$ and $9$.
\end{observation}

\subsection{Even dimension.}

\bigskip

\[
\begin{tabular}{l|l|l|l|}
\hline
\multicolumn{4}{|c|}{$\dim \;C^{1}\frak{g}=6$} \\ \hline
\multicolumn{1}{|l|}{Table $1$} & $\frak{g}$ & $\dim Z\left( \frak{g}\right) 
$ & $\left[ \frac{\frak{g}}{Z\left( g\right) }\right] $ \\ \hline
& $\frak{g}_{2m}^{1}$ & \multicolumn{1}{|c|}{$2$} & \multicolumn{1}{|c|}{$%
\frak{l}_{6}^{2}$} \\ \cline{2-4}
& $\frak{g}_{2m}^{2}$ & \multicolumn{1}{|c|}{$3$} & \multicolumn{1}{|c|}{$%
\frak{l}_{6}^{5}$} \\ \cline{2-4}
& $\frak{g}_{2m}^{3}$ & \multicolumn{1}{|c|}{$3$} & \multicolumn{1}{|c|}{$%
\frak{l}_{5}^{1}$} \\ \cline{2-4}
& $\frak{g}_{2m}^{4}$ & \multicolumn{1}{|c|}{$3$} & \multicolumn{1}{|c|}{$%
\frak{l}_{5}^{2}$} \\ \cline{2-4}
\end{tabular}
\]

\bigskip

\[
\begin{tabular}{l|c|c|c|c|}
\hline
\multicolumn{5}{|c|}{$\dim \;C^{1}\frak{g}=5$} \\ \hline
\multicolumn{1}{|l|}{Table $2$} & $\frak{g}$ & $\dim Z\left( \frak{g}\right) 
$ & $\left[ \frac{\frak{g}}{Z\left( \frak{g}\right) }\right] $ & $\left[ 
\frak{n}\right] $ \\ \hline
& $\frak{g}_{2m}^{6}$ & $1$ & $\frak{l}_{6}^{2}$ & $\frak{n}_{9}^{7}$ \\ 
\cline{2-5}
& $\frak{g}_{2m}^{7,\alpha }$ & $1$ & $\frak{l}_{6}^{2}$ & $\frak{n}%
_{13,\alpha }^{7}$ \\ \cline{2-5}
& $\frak{g}_{2m}^{8}$ & $1$ & $\frak{l}_{6}^{2}$ & $\frak{n}_{15}^{7}$ \\ 
\cline{2-5}
& $\frak{g}_{2m}^{9}$ & $1$ & $\frak{l}_{6}^{2}$ & $\frak{n}_{14}^{7}$ \\ 
\cline{2-5}
& $\frak{g}_{2m}^{10}$ & $1$ & $\frak{l}_{7}^{2}$ & $\frak{n}_{17}^{7}$ \\ 
\cline{2-5}
& $\frak{g}_{2m}^{11}$ & $1$ & $\frak{l}_{7}^{2}$ & $\frak{n}_{16}^{7}$ \\ 
\cline{2-5}
& $\frak{g}_{2m}^{12}$ & $2$ & $\frak{l}_{6}^{5}$ & $\frak{n}_{21}^{7}$ \\ 
\cline{2-5}
& $\frak{g}_{2m}^{13}$ & $2$ & $\frak{l}_{6}^{5}$ & $\frak{n}_{19}^{7}$ \\ 
\cline{2-5}
& $\frak{g}_{2m}^{14}$ & $2$ & $\frak{l}_{6}^{0}$ & $\frak{n}_{18}^{7}$ \\ 
\cline{2-5}
& $\frak{g}_{2m}^{15}$ & $2$ & $\frak{l}_{6}^{0}$ & $\frak{n}_{20}^{7}$ \\ 
\cline{2-5}
& $\frak{g}_{2m}^{16}$ & $2$ & $\frak{l}_{5}^{1}$ & $\frak{n}_{21}^{7}$ \\ 
\cline{2-5}
& $\frak{g}_{2m}^{17}$ & $2$ & $\frak{l}_{5}^{1}$ & $\frak{n}_{19}^{7}$ \\ 
\cline{2-5}
& $\frak{g}_{2m}^{18}$ & $2$ & $\frak{l}_{5}^{2}$ & $\frak{n}_{18}^{7}$ \\ 
\cline{2-5}
& $\frak{g}_{2m}^{19}$ & $2$ & $\frak{l}_{5}^{2}$ & $\frak{n}_{20}^{7}$ \\ 
\cline{2-5}
& $\frak{g}_{2m}^{20}$ & $2$ & $\frak{l}_{6}^{5}$ & $\frak{n}_{21}^{7}$ \\ 
\cline{2-5}
& $\frak{g}_{2m}^{21}$ & $2$ & $\frak{l}_{6}^{5}$ & $\frak{n}_{19}^{7}$ \\ 
\cline{2-5}
& $\frak{g}_{2m}^{22}$ & $2$ & $\frak{l}_{6}^{0}$ & $\frak{n}_{18}^{7}$ \\ 
\cline{2-5}
& $\frak{g}_{2m}^{23}$ & $2$ & $\frak{l}_{6}^{0}$ & $\frak{n}_{20}^{7}$ \\ 
\cline{2-5}
\end{tabular}
\]

\begin{remark}
As a consequence of the classification in dimension $7$ ( see [3]) we have $%
\frak{g}_{2m}^{7,\alpha }\simeq \frak{g}_{2m}^{7,\alpha ^{\prime }}$ if and
only if $\alpha ^{\prime }=\pm \alpha .$\ \newline
The pairs $\left\{ \left( 12,20\right) ,\left( 13,21\right) ,\left(
14,22\right) ,\left( 15,23\right) \right\} $ give the same entries. So they
have to be distinguished otherwise. See table $8$ for its distinction.
\end{remark}

\bigskip
\[
\begin{tabular}{c|c|c|c|c}
\hline
\multicolumn{5}{|c}{$\dim \;C^{1}\frak{g}=4$ and $C^{1}\frak{g}$ is not
abelian} \\ \hline
\multicolumn{1}{|c|}{Table $3$\ } & $\frak{g}$ & $\dim Z\left( \frak{g}%
\right) $ & $\left[ \frac{\frak{g}}{Z\left( \frak{g}\right) }\right] $ & 
\multicolumn{1}{|c|}{$\dim Der\left( \frak{g}\right) $} \\ \hline
& $\frak{g}_{2m}^{24}$ & $1$ & $\frak{l}_{6}^{5}$ & \multicolumn{1}{|c|}{$%
2m^{2}-9m+19$} \\ \cline{2-5}
& $\frak{g}_{2m}^{25}$ & $1$ & $\frak{l}_{6}^{0}$ & \multicolumn{1}{|c|}{$%
2m^{2}-9m+17$} \\ \cline{2-5}
& $\frak{g}_{2m}^{26}$ & $1$ & $\frak{l}_{6}^{5}$ & \multicolumn{1}{|c|}{$%
2m^{2}-9m+18$} \\ \cline{2-5}
& $\frak{g}_{2m}^{27}$ & $1$ & $\frak{l}_{6}^{0}$ & \multicolumn{1}{|c|}{$%
2m^{2}-9m+16$} \\ \cline{2-5}
& $\frak{g}_{2m}^{28}$ & $1$ & $\frak{l}_{5}^{1}$ &  \\ \cline{2-4}
& $\frak{g}_{2m}^{29}$ & $1$ & $\frak{l}_{5}^{2}$ &  \\ \cline{2-4}
\end{tabular}
\]

\bigskip
\textbf{\ } 
\[
\begin{tabular}{l|l|l|l|l}
\hline
\multicolumn{5}{|c}{$\dim \;C^{1}\frak{g}=4$ and $C^{1}\frak{g}$ \ is abelian
} \\ \hline
\multicolumn{1}{|l|}{Table $4$} & $\frak{g}$ & $\dim Z\left( \frak{g}\right) 
$ & $\left[ \frac{\frak{g}}{Z\left( \frak{g}\right) }\right] $ & 
\multicolumn{1}{|l|}{$\dim Der\left( \frak{g}\right) $} \\ \hline
& $\frak{g}_{2m}^{30}$ & $1$ & \multicolumn{1}{|c|}{$\frak{l}_{7}^{3,0}$} & 
\multicolumn{1}{|l|}{$2m^{2}-11m+27$} \\ \cline{2-5}
& $\frak{g}_{2m}^{31}$ & $1$ & \multicolumn{1}{|c|}{$\frak{l}_{7}^{3,0}$} & 
\multicolumn{1}{|l|}{$2m^{2}-11m+26$} \\ \cline{2-5}
& $\frak{g}_{2m}^{32}$ & $1$ & \multicolumn{1}{|c|}{$\frak{l}_{7}^{3,0}$} & 
\multicolumn{1}{|l|}{$2m^{2}-11m+26$} \\ \cline{2-5}
& $\frak{g}_{2m}^{33}$ & $1$ & \multicolumn{1}{|c|}{$\frak{l}_{7}^{3,0}$} & 
\multicolumn{1}{|l|}{$2m^{2}-11m+25$} \\ \cline{2-5}
& $\frak{g}_{2m}^{34}$ & $1$ & \multicolumn{1}{|c|}{$\frak{l}_{7}^{3,0}$} & 
\multicolumn{1}{|l|}{$2m^{2}-11m+27$} \\ \cline{2-5}
& $\frak{g}_{2m}^{35}$ & $1$ & \multicolumn{1}{|c|}{$\frak{l}_{7}^{3,0}$} & 
\multicolumn{1}{|l|}{$2m^{2}-11m+25$} \\ \cline{2-5}
& $\frak{g}_{2m}^{36}$ & $1$ & \multicolumn{1}{|c|}{$\frak{l}_{7}^{3,0}$} & 
\multicolumn{1}{|l|}{$2m^{2}-11m+27$} \\ \cline{2-5}
& $\frak{g}_{2m}^{37}$ & $1$ & \multicolumn{1}{|c|}{$\frak{l}_{7}^{3,0}$} & 
\multicolumn{1}{|l|}{$2m^{2}-11m+24$} \\ \cline{2-5}
& $\frak{g}_{2m}^{38}$ & $1$ & \multicolumn{1}{|c|}{$\frak{l}_{6}^{4}$} & 
\multicolumn{1}{|l|}{$2m^{2}-9m+18$} \\ \cline{2-5}
& $\frak{g}_{2m}^{39}$ & $1$ & \multicolumn{1}{|c|}{$\frak{l}_{6}^{4}$} & 
\multicolumn{1}{|l|}{$2m^{2}-9m+17$} \\ \cline{2-5}
& $\frak{g}_{2m}^{40}$ & $1$ & \multicolumn{1}{|c|}{$\frak{l}_{6}^{4}$} & 
\multicolumn{1}{|l|}{$2m^{2}-9m+19$} \\ \cline{2-5}
& $\frak{g}_{2m}^{41}$ & $1$ & \multicolumn{1}{|c|}{$\frak{l}_{6}^{4}$} & 
\multicolumn{1}{|l|}{$2m^{2}-9m+18$} \\ \cline{2-5}
& $\frak{g}_{2m}^{42}$ & $1$ & \multicolumn{1}{|c|}{$\frak{l}_{7}^{3,0}$} & 
\multicolumn{1}{|l|}{$2m^{2}-11m+25$} \\ \cline{2-5}
& $\frak{g}_{2m}^{43}$ & $1$ & \multicolumn{1}{|c|}{$\frak{l}_{6}^{4}$} & 
\multicolumn{1}{|l|}{$2m^{2}-9m+19$} \\ \cline{2-5}
& $\frak{g}_{2m}^{44}$ & $1$ & \multicolumn{1}{|c|}{$\frak{l}_{6}^{4}$} & 
\multicolumn{1}{|l|}{$2m^{2}-9m+18$} \\ \cline{2-5}
& $\frak{g}_{2m}^{45}$ & $1$ & \multicolumn{1}{|c|}{$\frak{l}_{6}^{5}$} & 
\multicolumn{1}{|l|}{$2m^{2}-9m+20$} \\ \cline{2-5}
& $\frak{g}_{2m}^{46}$ & $1$ & \multicolumn{1}{|c|}{$\frak{l}_{6}^{0}$} & 
\multicolumn{1}{|l|}{$2m^{2}-9m+19$} \\ \cline{2-5}
& $\frak{g}_{2m}^{47}$ & $1$ & \multicolumn{1}{|c|}{$\frak{l}_{6}^{5}$} & 
\multicolumn{1}{|l|}{$2m^{2}-9m+19$} \\ \cline{2-5}
& $\frak{g}_{2m}^{48}$ & $1$ & \multicolumn{1}{|c|}{$\frak{l}_{6}^{0}$} & 
\multicolumn{1}{|l|}{$2m^{2}-9m+18$} \\ \cline{2-5}
& $\frak{g}_{2m}^{49}$ & $1$ & \multicolumn{1}{|c|}{$\frak{l}_{6}^{5}$} & 
\multicolumn{1}{|l|}{$2m^{2}-9m+20$} \\ \cline{2-5}
& $\frak{g}_{2m}^{50}$ & $1$ & \multicolumn{1}{|c|}{$\frak{l}_{6}^{0}$} & 
\multicolumn{1}{|l|}{$2m^{2}-9m+18$} \\ \cline{2-5}\cline{5-5}
& $\frak{g}_{2m}^{51}$ & $1$ & \multicolumn{1}{|c|}{$\frak{l}_{5}^{1}$} & 
\\ \cline{2-4}
& $\frak{g}_{2m}^{52}$ & $1$ & \multicolumn{1}{|c|}{$\frak{l}_{5}^{2}$} & 
\\ \cline{2-4}
& $\frak{g}_{2m}^{53}$ & $1$ & \multicolumn{1}{|c|}{$\frak{l}_{5}^{2}$} & 
\\ \cline{2-4}
\end{tabular}
\]

\begin{remark}
The algebras $\frak{g}^{51},\frak{g}^{52},\frak{g}^{53}$ are distinguished
by their unique class of filiform ideals in dimension six, which are
respectively $\frak{N}_{6,5}$ ,$\frak{N}_{6,4}$ ,$\frak{N}_{6,3}.$ For this
reason it was not necessary to calculate the dimension of its algebra of
derivations. For the algebras non distinguished by the table above see table
9.
\end{remark}

\subsection{Odd dimension}

\bigskip

\[
\begin{tabular}{c|c|c|c|c|}
\hline
\multicolumn{5}{|c|}{$\dim \;C^{1}\frak{g}=5$} \\ \hline
\multicolumn{1}{|c|}{Table $5$} & $\frak{g}$ & $\dim Z\left( \frak{g}\right) 
$ & $\left[ \frac{\frak{g}}{Z\left( \frak{g}\right) }\right] $ & $\left[ 
\frak{n}\right] $ \\ \hline
& $\frak{g}_{2m+1}^{55}$ & $2$ & $\frak{l}_{6}^{5}$ & $\frak{n}_{12}^{7}$ \\ 
\cline{2-5}
& $\frak{g}_{2m+1}^{56}$ & $1$ & $\frak{l}_{7}^{2}$ & $\frak{n}_{17}^{7}$ \\ 
\cline{2-5}
& $\frak{g}_{2m+1}^{57}$ & $1$ & $\frak{l}_{7}^{2}$ & $\frak{n}_{16}^{7}$ \\ 
\cline{2-5}
& $\frak{g}_{2m+1}^{58}$ & $2$ & $\frak{l}_{6}^{5}$ & $\frak{n}_{21}^{7}$ \\ 
\cline{2-5}
& $\frak{g}_{2m+1}^{59}$ & $2$ & $\frak{l}_{7}^{2}$ & $\frak{n}_{19}^{7}$ \\ 
\cline{2-5}
& $\frak{g}_{2m+1}^{60}$ & $2$ & $\frak{l}_{6}^{0}$ & $\frak{n}_{18}^{7}$ \\ 
\cline{2-5}
& $\frak{g}_{2m+1}^{61}$ & $2$ & $\frak{l}_{6}^{0}$ & $\frak{n}_{20}^{7}$ \\ 
\cline{2-5}
& $\frak{g}_{2m+1}^{62}$ & $2$ & $\frak{l}_{5}^{2}$ & $\frak{n}_{10}^{7}$ \\ 
\cline{2-5}
& $\frak{g}_{2m+1}^{63}$ & $2$ & $\frak{l}_{5}^{1}$ & $\frak{n}_{10}^{7}$ \\ 
\cline{2-5}
& $\frak{g}_{2m+1}^{64}$ & $2$ & $\frak{l}_{5}^{1}$ & $\frak{n}_{11}^{7}$ \\ 
\cline{2-5}
& $\frak{g}_{2m+1}^{65}$ & $1$ & $\frak{l}_{6}^{2}$ & $\frak{n}_{9}^{7}$ \\ 
\cline{2-5}
& $\frak{g}_{2m+1}^{66,\alpha }$ & $1$ & $\frak{l}_{6}^{2}$ & $\frak{n}%
_{13,\alpha }^{7}$ \\ \cline{2-5}
& $\frak{g}_{2m+1}^{67}$ & $1$ & $\frak{l}_{6}^{2}$ & $\frak{n}_{15}^{7}$ \\ 
\cline{2-5}
& $\frak{g}_{2m+1}^{68}$ & $1$ & $\frak{l}_{6}^{2}$ & $\frak{n}_{14}^{7}$ \\ 
\cline{2-5}
& $\frak{g}_{2m+1}^{69}$ & $1$ & $\frak{l}_{6}^{3}$ & $\frak{n}_{17}^{7}$ \\ 
\cline{2-5}
& $\frak{g}_{2m+1}^{70}$ & $1$ & $\frak{l}_{6}^{3}$ & $\frak{n}_{16}^{7}$ \\ 
\cline{2-5}
& $\frak{g}_{2m+1}^{71}$ & $2$ & $\frak{l}_{5}^{1}$ & $\frak{n}_{21}^{7}$ \\ 
\cline{2-5}
& $\frak{g}_{2m+1}^{72}$ & $2$ & $\frak{l}_{5}^{1}$ & $\frak{n}_{19}^{7}$ \\ 
\cline{2-5}
& $\frak{g}_{2m+1}^{73}$ & $2$ & $\frak{l}_{5}^{2}$ & $\frak{n}_{18}^{7}$ \\ 
\cline{2-5}
& $\frak{g}_{2m+1}^{74}$ & $2$ & $\frak{l}_{5}^{2}$ & $\frak{n}_{20}^{7}$ \\ 
\cline{2-5}
& $\frak{g}_{2m+1}^{75}$ & $2$ & $\frak{l}_{6}^{5}$ & $\frak{n}_{21}^{7}$ \\ 
\cline{2-5}
& $\frak{g}_{2m+1}^{76}$ & $2$ & $\frak{l}_{7}^{2}$ & $\frak{n}_{19}^{7}$ \\ 
\cline{2-5}
& $\frak{g}_{2m+1}^{77}$ & $2$ & $\frak{l}_{6}^{0}$ & $\frak{n}_{18}^{7}$ \\ 
\cline{2-5}
& $\frak{g}_{2m+1}^{78}$ & $2$ & $\frak{l}_{6}^{0}$ & $\frak{n}_{20}^{7}$ \\ 
\cline{2-5}
\end{tabular}
\]

\begin{remark}
We observe that $\frak{g}_{2m+1}^{66,\alpha }\simeq \frak{g}%
_{2m+1}^{66,\alpha ^{\prime }}$ if and only if $\alpha ^{\prime }=\pm \alpha
.$ \newline
The pairs $\left\{ \left( 58,75\right) ,\left( 59,76\right) ,\left(
60,77\right) ,\left( 61,78\right) \right\} $ must be distinguished. See
Table 8.
\end{remark}

\bigskip

\[
\begin{tabular}{l|c|c|c|c}
\hline
\multicolumn{5}{|c}{$\dim \;C^{1}\frak{g}=4$ and $\;C^{1}\frak{g\;}$ is not
abelian} \\ \hline
\multicolumn{1}{|l|}{Table $6$} & $\frak{g}$ & $\dim Z\left( \frak{g}\right) 
$ & $\left[ \frac{\frak{g}}{Z\left( \frak{g}\right) }\right] $ & 
\multicolumn{1}{|c|}{$\dim Der\left( \frak{g}\right) $} \\ \hline
& $\frak{g}_{2m+1}^{79}$ & $1$ & $\frak{l}_{6}^{5}$ & \multicolumn{1}{|c|}{$%
2m^{2}-7m+16$} \\ \cline{2-5}
& $\frak{g}_{2m+1}^{80}$ & $1$ & $\frak{l}_{6}^{0}$ & \multicolumn{1}{|c|}{$%
2m^{2}-7m+13$} \\ \cline{2-5}
& $\frak{g}_{2m+1}^{81}$ & $1$ & $\frak{l}_{6}^{5}$ & \multicolumn{1}{|c|}{$%
2m^{2}-7m+13$} \\ \cline{2-5}
& $\frak{g}_{2m+1}^{82}$ & $1$ & $\frak{l}_{6}^{5}$ & \multicolumn{1}{|c|}{$%
2m^{2}-7m+14$} \\ \cline{2-5}
& $\frak{g}_{2m+1}^{83}$ & $1$ & $\frak{l}_{6}^{0}$ & \multicolumn{1}{|c|}{$%
2m^{2}-7m+13$} \\ \cline{2-5}\cline{5-5}
& $\frak{g}_{2m+1}^{84}$ & $1$ & $\frak{l}_{5}^{1}$ &  \\ \cline{2-4}
& $\frak{g}_{2m+1}^{85}$ & $1$ & $\frak{l}_{5}^{2}$ &  \\ \cline{2-4}
\end{tabular}
\]

\bigskip

\[
\begin{tabular}{l|c|c|c|c}
\hline
\multicolumn{5}{|c}{$\dim \;C^{1}\frak{g}=4$ and $C^{1}\frak{g}$ is abelian}
\\ \hline
\multicolumn{1}{|l|}{Table $7$} & $\frak{g}$ & $\dim Z\left( \frak{g}\right) 
$ & $\left[ \frac{\frak{g}}{Z\left( \frak{g}\right) }\right] $ & 
\multicolumn{1}{|c|}{$\dim Der\left( \frak{g}\right) $} \\ \hline
& $\frak{g}_{2m+1}^{86}$ & $1$ & $\frak{l}_{7}^{3,0}$ & \multicolumn{1}{|c|}{%
$2m^{2}-9m+22$} \\ \cline{2-5}
& $\frak{g}_{2m+1}^{87}$ & $1$ & $\frak{l}_{7}^{3,0}$ & \multicolumn{1}{|c|}{%
$2m^{2}-9m+19$} \\ \cline{2-5}
& $\frak{g}_{2m+1}^{88}$ & $1$ & $\frak{l}_{7}^{3,0}$ & \multicolumn{1}{|c|}{%
$2m^{2}-9m+22$} \\ \cline{2-5}
& $\frak{g}_{2m+1}^{89}$ & $1$ & $\frak{l}_{7}^{3,0}$ & \multicolumn{1}{|c|}{%
$2m^{2}-9m+21$} \\ \cline{2-5}
& $\frak{g}_{2m+1}^{90}$ & $1$ & $\frak{l}_{7}^{3,0}$ & \multicolumn{1}{|c|}{%
$2m^{2}-9m+20$} \\ \cline{2-5}
& $\frak{g}_{2m+1}^{91}$ & $1$ & $\frak{l}_{7}^{3,0}$ & \multicolumn{1}{|c|}{%
$2m^{2}-7m+13$} \\ \cline{2-5}
& $\frak{g}_{2m+1}^{92}$ & $1$ & $\frak{l}_{7}^{3,0}$ & \multicolumn{1}{|c|}{%
$2m^{2}-9m+20$} \\ \cline{2-5}
& $\frak{g}_{2m+1}^{93}$ & $1$ & $\frak{l}_{6}^{4}$ & \multicolumn{1}{|c|}{$%
2m^{2}-5m+8$} \\ \cline{2-5}
& $\frak{g}_{2m+1}^{94}$ & $1$ & $\frak{l}_{6}^{4}$ & \multicolumn{1}{|c|}{$%
2m^{2}-7m+15$} \\ \cline{2-5}
& $\frak{g}_{2m+1}^{95}$ & $1$ & $\frak{l}_{6}^{4}$ & \multicolumn{1}{|c|}{$%
2m^{2}-7m+13$} \\ \cline{2-5}
& $\frak{g}_{2m+1}^{96}$ & $1$ & $\frak{l}_{6}^{5}$ & \multicolumn{1}{|c|}{$%
2m^{2}-7m+16$} \\ \cline{2-5}
& $\frak{g}_{2m+1}^{97}$ & $1$ & $\frak{l}_{6}^{0}$ & \multicolumn{1}{|c|}{$%
2m^{2}-7m+15$} \\ \cline{2-5}
& $\frak{g}_{2m+1}^{98}$ & $1$ & $\frak{l}_{6}^{0}$ & \multicolumn{1}{|c|}{$%
2m^{2}-7m+14$} \\ \cline{2-5}
& $\frak{g}_{2m+1}^{99}$ & $1$ & $\frak{l}_{6}^{5}$ & \multicolumn{1}{|c|}{$%
2m^{2}-7m+16$} \\ \cline{2-5}
& $\frak{g}_{2m+1}^{100}$ & $1$ & $\frak{l}_{6}^{5}$ & \multicolumn{1}{|c|}{$%
2m^{2}-7m+15$} \\ \cline{2-5}
& $\frak{g}_{2m+1}^{101}$ & $1$ & $\frak{l}_{5}^{1}$ &  \\ \cline{2-4}
& $\frak{g}_{2m+1}^{102}$ & $1$ & $\frak{l}_{5}^{2}$ &  \\ \cline{2-4}
& $\frak{g}_{2m+1}^{103}$ & $1$ & $\frak{l}_{5}^{1}$ &  \\ \cline{2-4}
\end{tabular}
\]

\begin{remark}
The algebras $\frak{g}_{2m+1}^{101},\frak{g}_{2m+1}^{102}$ and $\frak{g}%
_{2m+1}^{103}$ are distinguished by their unique class of filiform ideals in
dimension six, which are respectively $\frak{N}_{6,5}$ ,$\frak{N}_{6,3}$ ,$%
\frak{N}_{6,4}$ .
\end{remark}

$
\begin{tabular}{|c|c}
\cline{1-1}
Table $8$ &  \\ \hline
$\frak{g}_{n}$ & \multicolumn{1}{|c|}{$p_{c}\left( \lambda \right) $} \\ 
\hline
$\frak{g}_{2m}^{12}$ & \multicolumn{1}{|l|}{$\left( f_{1}^{1}-\lambda
\right) \left( f_{2}^{2}-\lambda \right) \left( f_{1}^{1}+f_{2}^{2}-\lambda
\right) ^{2}\left( f_{1}^{1}+2f_{2}^{2}-\lambda \right) \left(
2f_{1}^{1}+f_{2}^{2}-\lambda \right) \left( 3f_{1}^{1}+f_{2}^{2}-\lambda
\right) \left( 4f_{1}^{1}+f_{2}^{2}-\lambda \right) \left(
3f_{1}^{1}-\lambda \right) \left( 4f_{1}^{1}-\lambda \right) \prod_{i\geq
5}\left( \mu _{i}-\lambda \right) $} \\ \hline
$\frak{g}_{2m}^{13}$ & \multicolumn{1}{|l|}{$\left( f_{1}^{1}-\lambda
\right) ^{2}\left( 2f_{1}^{1}-\lambda \right) ^{2}\left( 3f_{1}^{1}-\lambda
\right) ^{3}\left( 4f_{1}^{1}-\lambda \right) ^{2}\left( 5f_{1}^{1}-\lambda
\right) \prod_{i\geq 5}\left( \mu _{i}-\lambda \right) $} \\ \hline
$\frak{g}_{2m}^{14}$ & \multicolumn{1}{|l|}{$\left( f_{1}^{1}-\lambda
\right) \left( 2f_{1}^{1}-\lambda \right) \left( 3f_{1}^{1}-\lambda \right)
^{2}\left( 4f_{1}^{1}-\lambda \right) ^{2}\left( 5f_{1}^{1}-\lambda \right)
^{2}\left( 6f_{1}^{1}-\lambda \right) \left( 7f_{1}^{1}-\lambda \right)
\prod_{i\geq 5}\left( \mu _{i}-\lambda \right) $} \\ \hline
$\frak{g}_{2m}^{15}$ & \multicolumn{1}{|l|}{$\left( f_{1}^{1}-\lambda
\right) \left( 2f_{1}^{1}-\lambda \right) \left( 3f_{1}^{1}-\lambda \right)
^{3}\left( 4f_{1}^{1}-\lambda \right) ^{2}\left( 5f_{1}^{1}-\lambda \right)
^{2}\left( 6f_{1}^{1}-\lambda \right) \prod_{i\geq 5}\left( \mu _{i}-\lambda
\right) $} \\ \hline
$\frak{g}_{2m}^{20}$ & \multicolumn{1}{|l|}{$\left( f_{1}^{1}-\lambda
\right) \left( f_{2}^{2}-\lambda \right) \left( f_{1}^{1}+f_{2}^{2}-\lambda
\right) \left( 2f_{1}^{1}+f_{2}^{2}-\lambda \right) \left(
3f_{1}^{1}+f_{2}^{2}-\lambda \right) \left( 4f_{1}^{1}+f_{2}^{2}-\lambda
\right) \left( f_{1}^{1}+2f_{2}^{2}-\lambda \right) \left(
3f_{1}^{1}-\lambda \right) \prod_{i\geq 3}\left( \mu _{i}-\lambda \right) $}
\\ \hline
$\frak{g}_{2m}^{21}$ & \multicolumn{1}{|l|}{$\left( f_{1}^{1}-\lambda
\right) ^{2}\left( 2f_{1}^{1}-\lambda \right) \left( 3f_{1}^{1}-\lambda
\right) ^{3}\left( 4f_{1}^{1}-\lambda \right) \left( 5f_{1}^{1}-\lambda
\right) \prod_{i\geq 3}\left( \mu _{i}-\lambda \right) $} \\ \hline
$\frak{g}_{2m}^{22}$ & \multicolumn{1}{|l|}{$\left( f_{1}^{1}-\lambda
\right) \left( 2f_{1}^{1}-\lambda \right) \left( 3f_{1}^{1}-\lambda \right)
^{2}\left( 4f_{1}^{1}-\lambda \right) \left( 5f_{1}^{1}-\lambda \right)
^{2}\left( 7f_{1}^{1}-\lambda \right) \prod_{i\geq 3}\left( \mu _{i}-\lambda
\right) $} \\ \hline
$\frak{g}_{2m}^{23}$ & \multicolumn{1}{|l|}{$\left( f_{1}^{1}-\lambda
\right) \left( 2f_{1}^{1}-\lambda \right) \left( 3f_{1}^{1}-\lambda \right)
^{2}\left( 4f_{1}^{1}-\lambda \right) \left( 5f_{1}^{1}-\lambda \right)
^{2}\left( 6f_{1}^{1}-\lambda \right) \prod_{i\geq 3}\left( \mu _{i}-\lambda
\right) $} \\ \hline
$\frak{g}_{2m+1}^{58}$ & \multicolumn{1}{|l|}{$\left( f_{1}^{1}-\lambda
\right) \left( f_{2}^{2}-\lambda \right) \left( f_{1}^{1}+f_{2}^{2}-\lambda
\right) ^{2}\left( 2f_{1}^{1}+f_{2}^{2}-\lambda \right) \left(
3f_{1}^{1}+f_{2}^{2}-\lambda \right) \left( 4f_{1}^{1}+f_{2}^{2}-\lambda
\right) \left( f_{1}^{1}+2f_{2}^{2}-\lambda \right) \left(
3f_{1}^{1}-\lambda \right) \prod_{i\geq 4}\left( \mu _{i}-\lambda \right) $}
\\ \hline
$\frak{g}_{2m+1}^{59}$ & \multicolumn{1}{|l|}{$\left( f_{1}^{1}-\lambda
\right) \left( f_{2}^{2}-\lambda \right) \left( f_{1}^{1}+f_{2}^{2}-\lambda
\right) \left( 2f_{1}^{1}+f_{2}^{2}-\lambda \right) \left(
3f_{1}^{1}+f_{2}^{2}-\lambda \right) \left( 3f_{1}^{1}+2f_{2}^{2}-\lambda
\right) \left( f_{1}^{1}+2f_{2}^{2}-\lambda \right) \left(
2f_{1}^{1}-\lambda \right) \left( 3f_{1}^{1}-\lambda \right) \prod_{i\geq
4}\left( \mu _{i}-\lambda \right) $} \\ \hline
$\frak{g}_{2m+1}^{60}$ & \multicolumn{1}{|l|}{$\left( f_{1}^{1}-\lambda
\right) \left( 2f_{1}^{1}-\lambda \right) \left( 3f_{1}^{1}-\lambda \right)
^{2}\left( 4f_{1}^{1}-\lambda \right) ^{2}\left( 5f_{1}^{1}-\lambda \right)
^{2}\left( 6f_{1}^{1}-\lambda \right) \prod_{i\geq 4}\left( \mu _{i}-\lambda
\right) $} \\ \hline
$\frak{g}_{2m+1}^{61}$ & \multicolumn{1}{|l|}{$\left( f_{1}^{1}-\lambda
\right) \left( 2f_{1}^{1}-\lambda \right) \left( 3f_{1}^{1}-\lambda \right)
^{3}\left( 4f_{1}^{1}-\lambda \right) \left( 5f_{1}^{1}-\lambda \right)
^{2}\left( 6f_{1}^{1}-\lambda \right) \prod_{i\geq 4}\left( \mu _{i}-\lambda
\right) $} \\ \hline
$\frak{g}_{2m+1}^{75}$ & \multicolumn{1}{|l|}{$\left( f_{1}^{1}-\lambda
\right) \left( f_{2}^{2}-\lambda \right) \left( f_{1}^{1}+f_{2}^{2}-\lambda
\right) \left( f_{1}^{1}+2f_{2}^{2}-\lambda \right) \left(
2f_{1}^{1}+f_{2}^{2}-\lambda \right) \left( 3f_{1}^{1}+f_{2}^{2}-\lambda
\right) \left( 4f_{1}^{1}+f_{2}^{2}-\lambda \right) \left(
3f_{1}^{1}-\lambda \right) \left( 4f_{1}^{1}-\lambda \right) \prod_{i\geq
4}\left( \mu _{i}-\lambda \right) $} \\ \hline
$\frak{g}_{2m+1}^{76}$ & \multicolumn{1}{|l|}{$\left( f_{1}^{1}-\lambda
\right) \left( f_{2}^{2}-\lambda \right) \left( f_{1}^{1}+f_{2}^{2}-\lambda
\right) \left( 2f_{1}^{1}+f_{2}^{2}-\lambda \right) \left(
f_{1}^{1}+2f_{2}^{2}-\lambda \right) \left( 2f_{1}^{1}+2f_{2}^{2}-\lambda
\right) \left( 3f_{1}^{1}+f_{2}^{2}-\lambda \right) \left(
3f_{1}^{1}+2f_{2}^{2}-\lambda \right) \left( 3f_{1}^{1}-\lambda \right)
\prod_{i\geq 4}\left( \mu _{i}-\lambda \right) $} \\ \hline
$\frak{g}_{2m+1}^{77}$ & \multicolumn{1}{|l|}{$\left( f_{1}^{1}-\lambda
\right) \left( 2f_{1}^{1}-\lambda \right) \left( 3f_{1}^{1}-\lambda \right)
^{2}\left( 4f_{1}^{1}-\lambda \right) \left( 5f_{1}^{1}-\lambda \right)
^{2}\left( 6f_{1}^{1}-\lambda \right) \left( 7f_{1}^{1}-\lambda \right)
\prod_{i\geq 4}\left( \mu _{i}-\lambda \right) $} \\ \hline
$\frak{g}_{2m+1}^{78}$ & \multicolumn{1}{|l|}{$\left( f_{1}^{1}-\lambda
\right) \left( 2f_{1}^{1}-\lambda \right) \left( 3f_{1}^{1}-\lambda \right)
^{2}\left( 4f_{1}^{1}-\lambda \right) ^{2}\left( 5f_{1}^{1}-\lambda \right)
^{2}\left( 6f_{1}^{1}-\lambda \right) \prod_{i\geq 4}\left( \mu _{i}-\lambda
\right) $} \\ \hline
\end{tabular}
$

\bigskip 

$
\begin{tabular}{|c|c}
\cline{1-1}
Table $9$ &  \\ \hline
$\frak{g}$ & \multicolumn{1}{|c|}{$p_{c}\left( \lambda \right) $} \\ \hline
$\frak{g}_{2m}^{30}$ & \multicolumn{1}{|l|}{$\left( f_{1}^{1}-\lambda
\right) \left( f_{2}^{2}-\lambda \right) \left( f_{1}^{1}+f_{2}^{2}-\lambda
\right) \left( 2f_{1}^{1}+f_{2}^{2}-\lambda \right) ^{2}\left(
3f_{1}^{1}+f_{2}^{2}-\lambda \right) \left( 4f_{1}^{1}+f_{2}^{2}-\lambda
\right) \left( 2f_{1}^{1}-\lambda \right) \left( 3f_{1}^{1}-\lambda \right)
\left( 4f_{1}^{1}-\lambda \right) \prod_{i\geq 5}\left( \mu _{i}-\lambda
\right) $} \\ \hline
$\frak{g}_{2m}^{31}$ & \multicolumn{1}{|l|}{$\left( f_{1}^{1}-\lambda
\right) \left( 2f_{1}^{1}-\lambda \right) \left( 3f_{1}^{1}-\lambda \right)
\left( 4f_{1}^{1}-\lambda \right) \left( f_{2}^{2}-\lambda \right) \left(
f_{1}^{1}+f_{2}^{2}-\lambda \right) \left( 2f_{1}^{1}+f_{2}^{2}-\lambda
\right) \left( 3f_{1}^{1}+f_{2}^{2}-\lambda \right) ^{2}\left(
4f_{1}^{1}+f_{2}^{2}-\lambda \right) \prod_{i\geq 5}\left( \mu _{i}-\lambda
\right) $} \\ \hline
$\frak{g}_{2m}^{32}$ & \multicolumn{1}{|l|}{$\left( f_{1}^{1}-\lambda
\right) \left( 2f_{1}^{1}-\lambda \right) \left( 3f_{1}^{1}-\lambda \right)
^{2}\left( 4f_{1}^{1}-\lambda \right) ^{2}\left( 5f_{1}^{1}-\lambda \right)
^{2}\left( 6f_{1}^{1}-\lambda \right) \left( 7f_{1}^{1}-\lambda \right)
\prod_{i\geq 5}\left( \mu _{i}-\lambda \right) $} \\ \hline
$\frak{g}_{2m}^{33}$ & \multicolumn{1}{|l|}{$\left( f_{1}^{1}-\lambda
\right) \left( 2f_{1}^{1}-\lambda \right) ^{2}\left( 3f_{1}^{1}-\lambda
\right) ^{3}\left( 4f_{1}^{1}-\lambda \right) ^{2}\left( 5f_{1}^{1}-\lambda
\right) \left( 6f_{1}^{1}-\lambda \right) \prod_{i\geq 5}\left( \mu
_{i}-\lambda \right) $} \\ \hline
$\frak{g}_{2m}^{34}$ & \multicolumn{1}{|l|}{$\left( f_{1}^{1}-\lambda
\right) \left( f_{2}^{2}-\lambda \right) \left( f_{1}^{1}+f_{2}^{2}-\lambda
\right) ^{2}\left( 2f_{1}^{1}+f_{2}^{2}-\lambda \right) \left(
3f_{1}^{1}+f_{2}^{2}-\lambda \right) \left( 4f_{1}^{1}+f_{2}^{2}-\lambda
\right) \left( 2f_{1}^{1}-\lambda \right) \left( 3f_{1}^{1}-\lambda \right)
\left( 4f_{1}^{1}-\lambda \right) \prod_{i\geq 5}\left( \mu _{i}-\lambda
\right) $} \\ \hline
$\frak{g}_{2m}^{35}$ & \multicolumn{1}{|l|}{$\left( f_{1}^{1}-\lambda
\right) \left( 2f_{1}^{1}-\lambda \right) ^{2}\left( 3f_{1}^{1}-\lambda
\right) ^{3}\left( 4f_{1}^{1}-\lambda \right) ^{2}\left( 5f_{1}^{1}-\lambda
\right) \left( 6f_{1}^{1}-\lambda \right) \prod_{i\geq 5}\left( \mu
_{i}-\lambda \right) $} \\ \hline
$\frak{g}_{2m}^{36}$ & \multicolumn{1}{|l|}{$\left( f_{1}^{1}-\lambda
\right) \left( f_{2}^{2}-\lambda \right) \left( f_{1}^{1}+f_{2}^{2}-\lambda
\right) ^{2}\left( 2f_{1}^{1}+f_{2}^{2}-\lambda \right) ^{2}\left(
3f_{1}^{1}+f_{2}^{2}-\lambda \right) \left( 4f_{1}^{1}+f_{2}^{2}-\lambda
\right) \left( 2f_{1}^{1}-\lambda \right) \left( 3f_{1}^{1}-\lambda \right)
\prod_{i\geq 5}\left( \mu _{i}-\lambda \right) $} \\ \hline
$\frak{g}_{2m}^{38}$ & \multicolumn{1}{|l|}{$\left( f_{1}^{1}-\lambda
\right) \left( 2f_{1}^{1}-\lambda \right) ^{2}\left( 3f_{1}^{1}-\lambda
\right) \left( 4f_{1}^{1}-\lambda \right) ^{2}\left( 5f_{1}^{1}-\lambda
\right) \left( 6f_{1}^{1}-\lambda \right) \prod_{i\geq 3}\left( \mu
_{i}-\lambda \right) $} \\ \hline
$\frak{g}_{2m}^{40}$ & \multicolumn{1}{|l|}{$\left( f_{1}^{1}-\lambda
\right) \left( f_{2}^{2}-\lambda \right) \left( f_{1}^{1}+f_{2}^{2}-\lambda
\right) \left( 2f_{1}^{1}+f_{2}^{2}-\lambda \right) ^{2}\left(
3f_{1}^{1}+f_{2}^{2}-\lambda \right) \left( 4f_{1}^{1}+f_{2}^{2}-\lambda
\right) \left( 2f_{1}^{1}-\lambda \right) \prod_{i\geq 3}\left( \mu
_{i}-\lambda \right) $} \\ \hline
$\frak{g}_{2m}^{41}$ & \multicolumn{1}{|l|}{$\left( f_{1}^{1}-\lambda
\right) \left( 2f_{1}^{1}-\lambda \right) \left( 3f_{1}^{1}-\lambda \right)
\left( 4f_{1}^{1}-\lambda \right) \left( 5f_{1}^{1}-\lambda \right)
^{2}\left( 6f_{1}^{1}-\lambda \right) \left( 7f_{1}^{1}-\lambda \right)
\prod_{i\geq 3}\left( \mu _{i}-\lambda \right) $} \\ \hline
$\frak{g}_{2m}^{42}$ & \multicolumn{1}{|l|}{$\left( f_{1}^{1}-\lambda
\right) \left( f_{2}^{2}-\lambda \right) \left( f_{1}^{1}+f_{2}^{2}-\lambda
\right) \left( 2f_{1}^{1}+f_{2}^{2}-\lambda \right) \left(
3f_{1}^{1}+f_{2}^{2}-\lambda \right) \left( 4f_{1}^{1}+f_{2}^{2}-\lambda
\right) \left( 2f_{1}^{1}-\lambda \right) \left( 3f_{1}^{1}-\lambda \right)
\prod_{i\geq 3}\left( \mu _{i}-\lambda \right) $} \\ \hline
$\frak{g}_{2m}^{43}$ & \multicolumn{1}{|l|}{$\left( f_{1}^{1}-\lambda
\right) \left( f_{2}^{2}-\lambda \right) \left( f_{1}^{1}+f_{2}^{2}-\lambda
\right) \left( 2f_{1}^{1}+f_{2}^{2}-\lambda \right) \left(
3f_{1}^{1}+f_{2}^{2}-\lambda \right) \left( 4f_{1}^{1}+f_{2}^{2}-\lambda
\right) \left( 2f_{1}^{1}-\lambda \right) \left( 4f_{1}^{1}-\lambda \right)
\prod \left( \mu _{i}-\lambda \right) $} \\ \hline
$\frak{g}_{2m}^{45}$ & \multicolumn{1}{|c|}{$\left( f_{1}^{1}-\lambda
\right) \left( f_{2}^{2}-\lambda \right) \left( f_{1}^{1}+f_{2}^{2}-\lambda
\right) \left( 2f_{1}^{1}+f_{2}^{2}-\lambda \right) \left(
3f_{1}^{1}+f_{2}^{2}-\lambda \right) \left( 4f_{1}^{1}+f_{2}^{2}-\lambda
\right) \left( 3f_{1}^{1}-\lambda \right) \left( 4f_{1}^{1}-\lambda \right)
\prod_{i\geq 3}\left( \mu _{i}-\lambda \right) $} \\ \hline
$\frak{g}_{2m}^{48}$ & \multicolumn{1}{|c|}{$\left( f_{1}^{1}-\lambda
\right) \left( 2f_{1}^{1}-\lambda \right) \left( 3f_{1}^{1}-\lambda \right)
^{2}\left( 4f_{1}^{1}-\lambda \right) ^{2}\left( 5f_{1}^{1}-\lambda \right)
\left( 6f_{1}^{1}-\lambda \right) \prod_{i\geq 3}\left( \mu _{i}-\lambda
\right) $} \\ \hline
$\frak{g}_{2m}^{49}$ & \multicolumn{1}{|c|}{$\left( f_{1}^{1}-\lambda
\right) \left( f_{2}^{2}-\lambda \right) \left( f_{1}^{1}+f_{2}^{2}-\lambda
\right) ^{2}\left( 2f_{1}^{1}+f_{2}^{2}-\lambda \right) \left(
3f_{1}^{1}+f_{2}^{2}-\lambda \right) \left( 4f_{1}^{1}+f_{2}^{2}-\lambda
\right) \prod_{i\geq 3}\left( \mu _{i}-\lambda \right) $} \\ \hline
$\frak{g}_{2m}^{50}$ & \multicolumn{1}{|c|}{$\left( f_{1}^{1}-\lambda
\right) \left( 2f_{1}^{1}-\lambda \right) \left( 3f_{1}^{1}-\lambda \right)
^{3}\left( 4f_{1}^{1}-\lambda \right) \left( 5f_{1}^{1}-\lambda \right)
\left( 6f_{1}^{1}-\lambda \right) \prod_{i\geq 3}\left( \mu _{i}-\lambda
\right) $} \\ \hline
$\frak{g}_{2m+1}^{80}$ & \multicolumn{1}{|c|}{$\lambda ^{9}\prod_{i\geq
4}\left( \mu _{i}-\lambda \right) $} \\ \hline
$\frak{g}_{2m+1}^{83}$ & \multicolumn{1}{|c|}{$\lambda ^{7}\prod_{i\geq
2}\left( \mu _{i}-\lambda \right) $} \\ \hline
$\frak{g}_{2m+1}^{86}$ & \multicolumn{1}{|l|}{$\left( f_{1}^{1}-\lambda
\right) \left( f_{2}^{2}-\lambda \right) \left( f_{1}^{1}+f_{2}^{2}-\lambda
\right) ^{2}\left( 2f_{1}^{1}+f_{2}^{2}-\lambda \right) ^{2}\left(
3f_{1}^{1}+f_{2}^{2}-\lambda \right) \left( 4f_{1}^{1}+f_{2}^{2}-\lambda
\right) \left( 2f_{1}^{1}-\lambda \right) \left( 3f_{1}^{1}-\lambda \right)
\left( 4f_{1}^{1}-\lambda \right) \prod_{i\geq 5}\left( \mu _{i}-\lambda
\right) $} \\ \hline
$\frak{g}_{2m+1}^{88}$ & \multicolumn{1}{|c|}{$\left( f_{1}^{1}-\lambda
\right) \left( f_{2}^{2}-\lambda \right) \left( f_{1}^{1}+f_{2}^{2}-\lambda
\right) \left( 2f_{1}^{1}+f_{2}^{2}-\lambda \right) \left(
3f_{1}^{1}+f_{2}^{2}-\lambda \right) \left( 4f_{1}^{1}+f_{2}^{2}-\lambda
\right) \left( 2f_{1}^{1}-\lambda \right) \left( 3f_{1}^{1}-\lambda \right)
\left( 4f_{1}^{1}-\lambda \right) \prod_{i\geq 4}\left( \mu _{i}-\lambda
\right) $} \\ \hline
$\frak{g}_{2m+1}^{90}$ & \multicolumn{1}{|c|}{$\left( f_{1}^{1}-\lambda
\right) \left( 2f_{1}^{1}-\lambda \right) ^{2}\left( 3f_{1}^{1}-\lambda
\right) ^{2}\left( 4f_{1}^{1}-\lambda \right) ^{2}\left( 5f_{1}^{1}-\lambda
\right) \left( 6f_{1}^{1}-\lambda \right) \prod_{i\geq 4}\left( \mu
_{i}-\lambda \right) $} \\ \hline
$\frak{g}_{2m+1}^{92}$ & \multicolumn{1}{|c|}{$\left( f_{1}^{1}-\lambda
\right) \left( f_{2}^{2}-\lambda \right) \left( f_{1}^{1}+f_{2}^{2}-\lambda
\right) ^{2}\left( 2f_{1}^{1}+f_{2}^{2}-\lambda \right) \left(
3f_{1}^{1}+f_{2}^{2}-\lambda \right) \left( 4f_{1}^{1}+f_{2}^{2}-\lambda
\right) \left( 2f_{1}^{1}-\lambda \right) \left( 3f_{1}^{1}-\lambda \right)
\prod_{i\geq 4}\left( \mu _{i}-\lambda \right) $} \\ \hline
$\frak{g}_{2m+1}^{96}$ & \multicolumn{1}{|c|}{$\left( f_{1}^{1}-\lambda
\right) \left( f_{2}^{2}-\lambda \right) \left( f_{1}^{1}+f_{2}^{2}-\lambda
\right) ^{2}\left( 2f_{1}^{1}+f_{2}^{2}-\lambda \right) \left(
3f_{1}^{1}+f_{2}^{2}-\lambda \right) \left( 4f_{1}^{1}+f_{2}^{2}-\lambda
\right) \left( 3f_{1}^{1}-\lambda \right) \left( 4f_{1}^{1}-\lambda \right)
\prod_{i\geq 4}\left( \mu _{i}-\lambda \right) $} \\ \hline
$\frak{g}_{2m+1}^{99}$ & \multicolumn{1}{|c|}{$\left( f_{1}^{1}-\lambda
\right) \left( f_{2}^{2}-\lambda \right) \left( f_{1}^{1}+f_{2}^{2}-\lambda
\right) \left( 2f_{1}^{1}+f_{2}^{2}-\lambda \right) \left(
3f_{1}^{1}+f_{2}^{2}-\lambda \right) \left( 4f_{1}^{1}+f_{2}^{2}-\lambda
\right) \left( 3f_{1}^{1}-\lambda \right) \prod_{i\geq 2}\left( \mu
_{i}-\lambda \right) $} \\ \hline
\end{tabular}
$

\newpage
\section{\textbf{Proofs}}

\subsection{Preliminar results.}

We now prove the classification theorem. For simplicity and shortness we
omit several calculations.

\begin{lemma}
Let $\frak{g}_{n}$ be a $n$-dimensional complex $\left( n-5\right) $%
-filiform Lie algebra. Then the law $\mu _{n}$ is expressible as\newline
$\mu _{n}\left( X_{1},X_{j}\right) =X_{j+1},\;j=2,3,4,5.$\newline
$\mu _{n}\left( Y_{i},X_{4}\right) =d_{i}X_{6},\;i\in \{1,..,n-6\}$\newline
$\mu _{n}\left( Y_{i},X_{3}\right) =d_{i}X_{5}+e_{i}X_{6},\;i\in
\{1,..,n-6\} $\newline
$\mu _{n}\left( Y_{i},X_{2}\right) =d_{i}X_{4}+e_{i}X_{5}+f_{i}X_{6},\;i\in
\{1,..,n-6\}$\newline
$\mu _{n}\left( X_{5},X_{2}\right) =\mu _{n}\left( X_{3},X_{4}\right)
=\sum_{i=1}^{n-6}m^{k}Y_{k}+n_{6}X_{6};\;\mu _{n}\left( X_{4},X_{2}\right)
=o_{6}X_{6};$\newline
$\mu _{n}\left( X_{3},X_{2}\right)
=\sum_{i=1}^{n-6}p^{k}Y_{k}+o_{6}X_{5}+p_{6}X_{6};\;\mu _{n}\left(
Y_{i},Y_{j}\right) =a_{ij}X_{6},\;1\leq i,j\leq n-6$\newline
where the undefined brackets are zero or obtained by antisymmetry.%
\endproof%
\end{lemma}

For later use we introduce the following linear changes of basis:

\begin{itemize}
\item  Type I\newline
\[
\left\{ 
\begin{array}{l}
g_{1}\left( X_{2}\right)
=\sum_{i=2}^{6}a_{i}X_{i}+a_{7}Y_{1}+a_{2}Y_{2},\;a_{2}\neq 0 \\ 
g_{1}\left( Y_{1}\right) =b_{1}Y_{1}+b_{2}Y_{2}+b_{3}X_{6},\;b_{1}\neq
0\;\;\;\; \\ 
g_{1}\left( Y_{2}\right) =b_{4}Y_{2}+b_{5}X_{6},\;b_{4}\neq 0
\end{array}
\right. 
\]
where the undefined images are fixed by $g_{1}$ or obtained by the operator $%
ad\left( X_{2}\right) .$

\item  Type II 
\[
\left\{ 
\begin{array}{l}
g_{2}\left( X_{1}\right)
=\sum_{i=1}^{6}a_{i}X_{i}+a_{7}Y_{1}+a_{8}Y_{2},\;a_{1}\neq 0 \\ 
g_{2}\left( Y_{1}\right)
=b_{2}Y_{1}+b_{3}X_{6}+b_{4}X_{5}+b_{5}X_{4}+b_{6}X_{3}+b_{7}Y_{2} \\ 
g_{2}\left( Y_{2}\right) =c_{2}Y_{2}+c_{3}X_{6}+c_{4}X_{5}+c_{5}X_{4} \\ 
g_{2}\left( Y_{j}\right) =Y_{j}+\delta _{j}X_{5},\;j\geq 3
\end{array}
\right. 
\]
where the undefined images are fixed by $g_{2}$ or obtained by the operator $%
ad\left( X_{1}\right) .$

\item  Type III 
\[
\left\{ 
\begin{array}{l}
g_{3}\left( Y_{1}\right) =Y_{1}+\alpha _{1}Y_{4}+\alpha _{2}Y_{5} \\ 
g_{3}\left( Y_{2}\right) =Y_{2}+\alpha _{3}Y_{4}+\alpha _{4}Y_{5} \\ 
g_{3}\left( Y_{3}\right) =Y_{3}+\alpha _{5}Y_{4}+\alpha _{6}Y_{5} \\ 
g_{3}\left( Y_{4}\right) =Y_{4}+\alpha _{7}Y_{5}
\end{array}
\right. 
\]
where the undefined images are fixed by $g_{3}$ or obtained by the adjoint
operators$.$

\item  Type IV 
\[
\left\{ 
\begin{array}{l}
g_{4}\left( X_{2}\right) =X_{2}+\alpha Y_{3} \\ 
g_{4}\left( Y_{1}\right) =Y_{1}+\beta Y_{3} \\ 
g_{4}\left( Y_{2}\right) =Y_{2}+\gamma Y_{3}
\end{array}
\right. 
\]
where the undefined images are fixed by $g_{4}$ or obtained by adjoint
operators.
\end{itemize}

\begin{remark}
The new structure constants for the changes of type III and IV are easily
obtained. However, the first two types have to be commented.
\end{remark}

\begin{lemma}
For a change of basis of type I the new structure constants are 
\begin{eqnarray*}
\stackrel{\sim }{o_{6}} &=&a_{2}o_{6}-a_{7}d_{1} \\
\;\stackrel{\sim }{p_{6}} &=&a_{2}p_{6}-a_{7}e_{1}-a_{8}e_{2} \\
\stackrel{\sim }{d_{1}} &=&b_{1}d_{1} \\
\stackrel{\sim }{e_{1}} &=&b_{1}e_{1}+b_{2}e_{2} \\
\stackrel{\sim }{f_{1}} &=&b_{1}f_{1}+b_{2}f_{2}+\frac{b_{1}a_{8}a_{12}}{%
a_{2}} \\
\stackrel{\sim }{e_{2}} &=&b_{4}e_{2} \\
\stackrel{\sim }{f_{2}} &=&b_{4}f_{2}-\frac{a_{7}b_{4}a_{12}}{a_{2}} \\
\stackrel{\sim }{a_{12}} &=&b_{1}b_{4}a_{12} \\
\stackrel{\sim }{a_{1j}} &=&\frac{b_{1}a_{1j}+b_{2}a_{2j}}{a_{2}},\;j\geq 3\;
\\
\stackrel{\sim }{a_{2j}} &=&\frac{b_{4}a_{2j}}{a_{2}},\;j\geq
3\;\;\;\;\;\;\;\;\;\;%
\endproof%
\end{eqnarray*}
\end{lemma}

\begin{lemma}
For a change of basis of type II the new structure constants are 
\begin{eqnarray*}
\stackrel{\approx }{o_{6}} &=&\frac{o_{6}}{a_{1}^{2}} \\
\stackrel{\approx }{p_{6}} &=&\frac{%
a_{1}p_{6}+2a_{2}o_{6}^{2}-2a_{7}d_{1}o_{6}}{a_{1}^{4}} \\
\stackrel{\approx }{d_{1}} &=&\frac{b_{2}d_{1}}{a_{1}^{2}} \\
\stackrel{\approx }{e_{1}} &=&\frac{b_{2}e_{1}+b_{7}e_{2}}{a_{1}^{3}}+\frac{%
2b_{2}d_{1}}{a_{1}^{4}}\left( a_{2}o_{6}-a_{7}d_{1}\right) \\
\stackrel{\approx }{f_{1}} &=&\frac{5a_{7}^{2}b_{2}}{a_{1}^{6}}d_{1}^{3}-%
\frac{10a_{2}a_{7}o_{6}b_{2}}{a_{1}^{6}}d_{1}^{2}+\frac{b_{2}f_{1}+b_{7}f_{2}%
}{a_{1}^{4}}+\frac{3a_{2}o_{6}}{a_{1}^{5}}\left(
b_{2}e_{1}+b_{7}e_{2}\right) + \\
&&+\left( \frac{5a_{2}^{2}o_{6}^{2}b_{2}}{a_{1}^{6}}+\frac{2a_{2}b_{2}p_{6}}{%
a_{1}^{5}}-\frac{5a_{7}b_{2}e_{1}}{a_{1}^{5}}-\frac{3a_{7}b_{7}e_{2}}{%
a_{1}^{5}}-\frac{2b_{2}a_{8}e_{2}}{a_{1}^{5}}\right) d_{1} \\
\stackrel{\approx }{e_{2}} &=&\frac{c_{2}e_{2}}{a_{1}^{3}} \\
\stackrel{\approx }{f_{2}} &=&c_{2}\left( \frac{f_{2}}{a_{1}^{4}}+\frac{%
3a_{2}e_{2}o_{6}}{a_{1}^{5}}-\frac{3e_{2}a_{7}d_{1}}{a_{1}^{5}}\right) \\
\stackrel{\approx }{a_{12}} &=&\frac{b_{2}c_{2}a_{12}}{a_{1}^{4}} \\
\stackrel{\approx }{a_{1j}} &=&\frac{b_{2}a_{1j}+b_{7}a_{2j}}{a_{1}^{4}}-%
\frac{a_{2}b_{2}d_{1}e_{j}}{a_{1}^{5}},\;j\geq 3 \\
\stackrel{\approx }{a_{2j}} &=&\frac{c_{2}a_{2j}}{a_{1}^{4}}%
\;\;\;\;\;\;\;\;\;\;\;\;\;\;\;\;\;\;\;\;\;\;\;\;\;\;\;\;\;\;\;\;\;\;\;\;\;\;%
\;\;\;\;\;\;\;\;\;\;\;\;\;\;\;\;\;\;\;\;\;\;\;\;\;\;\;\;\;%
\endproof%
\end{eqnarray*}
\end{lemma}

\begin{remark}
From now on, we abbreviate the linear changes of basis with l.c.b.
\end{remark}

\subsection{Proof of the classification theorem}

\proof%
From now on, when we say algebra we mean an $\left( n-5\right) $-filiform
nonsplit Lie algebra.\newline
We begin considering the structure constants referring to the $Y_{i}$ in the
bracket $[X_{5},X_{2}]=[X_{3},X_{4}].$ This choice allows to separate the
cases by the dimension of the corresponding $C^{1}\frak{g}$ $.$ \newline

\subsubsection{There exists\textbf{\ }$m^{k}\neq 0.$}

We can suppose $m^{1}=1,m^{i}=0$ $\forall i\geq 2$ and $n_{6}=0.$ The Jacobi
conditions imply $d_{1}=e_{1}=a_{1j}=0;\;f_{1}=\sum_{k\geq
2}p^{k}d_{k};\;\sum p^{k}a_{ik}=0$ for $i,k\geq 2.$

\begin{enumerate}
\item  If $\exists p^{k}\neq 0$ with $k\neq 1$ we can suppose $%
p^{2}=1,p^{k}=0,\,\forall k\neq 2,$ and $p_{6}=0$ through a linear change of
basis. From the Jacobi conditions above we deduce $f_{1}=a_{i2}=0.$
Considering the characteristic vector $X_{1}^{\prime }=X_{1}+\alpha
X_{2},\;\alpha \neq 0,$ the characteristic sequence implies the system 
\[
\left\{ 
\begin{array}{c}
\;\;\;f_{j}=0,\;j\geq 2 \\ 
e_{j}o_{6}=0,\;j\geq 2
\end{array}
\right. 
\]
We distinguish the two cases

\begin{enumerate}
\item  If $e_{2}\neq 0$ we suppose $e_{2}=1$ and $e_{j}=0,\;\forall j\neq 2$
through the change $Y_{j}^{\prime }=Y_{j}-e_{j}Y_{2},\;j\neq 2.$ Reordering
the $Y_{i}$ we can suppose $[Y_{2t-1},Y_{2t}]=X_{6}$ for $2\leq t\leq \frac{%
n-6}{2};[Y_{i},Y_{j}]=0$ for the remaining. We obtain a unique class of
nonsplit Lie algebras in even dimension and isomorphic to $\frak{g}_{2m}^{1}.
$

\item  If $e_{2}=0$

\begin{enumerate}
\item  If there exists $e_{i}\neq 0$ with $i\geq 3$ we can suppose $e_{3}=1$
and $e_{j}=0$ for $j\neq 3.$ Reordering the $Y_{i}$ we obtain one algebra in
even and one algebra in odd dimension, which are respectively isomorphic to $%
\frak{g}_{2m}^{4}$ and $\frak{g}_{2m+1}^{54}.$

\item  If $e_{i}=0$ for $i\geq 3$ we obtain in an analogous way two even
dimensional algebras respectively isomorphic to $\frak{g}_{2m}^{2}$ and $%
\frak{g}_{2m}^{3}.$
\end{enumerate}

\begin{remark}
These are the Lie algebras with derived subalgebra of maximal dimension.
\end{remark}
\end{enumerate}

\item  $p^{k}=0$ for all $k.$ The characteristic sequence implies $%
f_{i}=e_{i}o_{6}=0$ for all $i.$ Moreover, the change of basis $%
X_{2}^{\prime }=X_{2}-\frac{p^{1}}{2}X_{4}$ allows to suppose $p^{1}=0.$

\begin{enumerate}
\item  If $\exists \,e_{i}\neq 0$ we suppose $e_{2}=1,e_{i}=0,\forall i\neq
2.$ A change of basis allows $p_{6}=0.$ There are two possibilities: an even
dimensional algebra isomorphic to $\frak{g}_{2m}^{5}$ and an even
dimensional one isomorphic to $\frak{g}_{2m+1}^{55}.$

\item  $e_{i}=0,\forall i.$

\begin{enumerate}
\item  If $o_{6}\neq 0$ we put $o_{6}=1$ and $p_{6}=0$ with a linear change
of basis. We obtain a unique algebra in odd dimension isomorphic to $\frak{g}%
_{2m+1}^{62}.$

\item  If $o_{6}=0$ there are two possibilities, depending on $p_{6}$ zero
or not. We obtain two odd dimensional algebras respectively isomorphic to $%
\frak{g}_{2m+1}^{63}$ and $\frak{g}_{2m+1}^{64}.$
\end{enumerate}
\end{enumerate}
\end{enumerate}

\subsubsection{$m^{k}=0,\;\forall k.$}

\begin{enumerate}
\item  If $p^{k}\neq 0$ for $k\geq 1$ we can suppose $p^{1}=1$ and $p^{i}=0$
for all $i\geq 2.$ The Jacobi conditions imply 
\[
d_{1}=0,\;a_{i1}=2d_{i}n_{6},\;i\geq 2
\]
From the characteristic sequence we deduce $d_{i}=0$ for all $i,$ so $a_{i1}$
from the Jacobi conditions.

\begin{enumerate}
\item  If $e_{i}\neq 0$ for an index $i$ we can suppose $e_{1}=1$ and the
remaining zero. The combination of a change of Jordan basis and a change of
type  $X_{1}^{\prime }=X_{1}+a_{1}X_{2}+a_{2}X_{3}+a_{3}Y_{1}$  with $%
a_{i}\in \Bbb{C}$ allows to take $n_{6}=0.$

\begin{enumerate}
\item  If $f_{1}\neq 0$ then $f_{1}=1$ through a linear change of basis.

\begin{enumerate}
\item  If $f_{i}=0$, $\forall i\geq 2$ we reorder the $\left\{
Y_{2},..,Y_{n-6}\right\} $ such that $[Y_{2t},Y_{2t+1}]=X_{6}$ $\,$for $%
1\leq t\leq \frac{n-6}{2}$ and the remaining brackets zero. The decisive
structure constant is $o_{6}.$ If it is zero we obtain an odd dimensional
Lie algebra isomorphic to $\frak{g}_{2m+1}^{65}.$ If not, $o_{6}=\alpha $ is
a parameter that can not be normalized by any change of basis. So we obtain
an infinite family of odd dimensional Lie algebras isomorphic to the family $%
\frak{g}_{2m+1}^{66,\alpha }.$

\item  $\exists \,f_{i}\neq 0$, $i\geq 2.$ Without loss of generality we can
choose $f_{2}=0$ and the remaining zero for $i\geq 3.$ It is easy deduce $%
a_{2j}=0,\,\forall j.$ Reordering the $\{Y_{3},..,Y_{n-6}\}$ in an
Heisenberg manner we obtain an even dimensional Lie algebra and an infinite
family of even dimensional algebras, which are respectively isomorphic to $%
\frak{g}_{2m}^{6}$ and $\frak{g}_{2m}^{7,\alpha }.$
\end{enumerate}

\item  Take $f_{1}=0.$

\begin{enumerate}
\item  If there is an index $i\geq 2$ such that $f_{i}\neq 0$ we can suppose 
$f_{2}=1$ and the remaining zero. Reordering the $\left\{
Y_{3},..,Y_{n-6}\right\} $ in an Heisenberg manner we obtain two even
dimensional Lie algebras, respectively isomorphic to $\frak{g}_{2m}^{8}$ and 
$\frak{g}_{2m}^{9}.$

\item  $f_{i}=0,\;\forall i.$ A similar reordering of the $Y_{i}$ gives two
Lie algebras in odd dimension isomorphic to $\frak{g}_{2m+1}^{67}$ and $%
\frak{g}_{2m+1}^{68}.$
\end{enumerate}
\end{enumerate}

\item  Let $e_{1}=0.$

\begin{enumerate}
\item  $f_{1}=1,\;f_{i}=0$ for $i\geq 2.$ A l.c.b. allows to suppose $%
o_{6}=0.$

\begin{enumerate}
\item  If $e_{i}\neq 0$ for $i\geq 2$ choose $e_{2}=1$ and the remaining
zero.\newline
There are two possible cases, depending on $a_{23}:$ if it is nonzero, we
obtain two algebras in odd dimension isomorphic to $\frak{g}_{2m+1}^{56}$
and $\frak{g}_{2m+1}^{57},$ and if it is zero, we obtain two even
dimensional algebras respectively isomorphic to $\frak{g}_{2m}^{10}$ and $%
\frak{g}_{2m}^{11}.$

\item  If $e_{i}=0,\;\forall i\geq 2$ a reordering of $\left\{
Y_{2},..,Y_{n-6}\right\} $ gives two algebras isomorphic to $\frak{g}%
_{2m+1}^{69}$ and $\frak{g}_{2m+1}^{70}.$
\end{enumerate}

\item  $f_{1}=0.$

\begin{enumerate}
\item  $\;e_{2}\neq 0$ and $e_{i}=0,\forall i\geq 3.$ With a linear change
of basis we can suppose $f_{2}=0.$\newline
A-1) If $\exists \,i\geq 3$ with $f_{i}\neq 0$ take $f_{3}=1$ and $f_{i}=0$
for $i\geq 4.$ A linear change allows to suppose $a_{3j}=0$ for all $j.$ If $%
a_{2j}=0$ for all $j$ we obtain the algebras $\frak{g}_{2m+1}^{75},\frak{g}%
_{2m+1}^{76},\frak{g}_{2m+1}^{77},\frak{g}_{2m+1}^{78}.$ If not, reorder the 
$\left\{ Y_{4},..,Y_{n-6}\right\} $ such that $a_{24}=1.$ We obtain the
algebras $\frak{g}_{2m}^{20},\frak{g}_{2m}^{21},\frak{g}_{2m}^{22}$ and $%
\frak{g}_{2m}^{23}.$\newline
A-2) $f_{i}=0\;$ $\forall i.$ Again, if $a_{23}=0$ we obtain the Lie
algebras $\frak{g}_{2m}^{12},\frak{g}_{2m}^{13},\frak{g}_{2m}^{14}$ and $%
\frak{g}_{2m}^{15},$ and if $a_{23}\neq 0$ we obtain the algebras $\frak{g}%
_{2m+1}^{58},\frak{g}_{2m+1}^{59},\frak{g}_{2m+1}^{60}$ and $\frak{g}%
_{2m+1}^{61}.$

\item  $e_{i}=0\;\forall i.$\newline
B-1) If there is an $i\geq 2$ with $f_{i}\neq 0$ take $f_{2}=1.$ A l.c.b.
allows to suppose $a_{2j}=0.$ We obtain four Lie algebras in even dimension
isomorphic to $\frak{g}_{2m}^{16},\frak{g}_{2m}^{17},\frak{g}_{2m}^{18}$ and 
$\frak{g}_{2m}^{19}.$\newline
B-2) If $f_{i}=0$ for $i\geq 2,$ the nonsplittness forces $a_{2j}$ for an
index $j.$ Reordering the $Y_{i}$ adequately, we can suppose $a_{23}=1.$ We
obtain odd dimensional algebras isomorphic to $\frak{g}_{2m+1}^{71},\frak{g}%
_{2m+1}^{72},\frak{g}_{2m+1}^{73}$ and $\frak{g}_{2m+1}^{74}.$
\end{enumerate}
\end{enumerate}
\end{enumerate}

\item  $p^{k}=0,\;\forall k.$\newline
The only Jacobi condition implies 
\begin{equation}
d_{i}n_{6}=0,\;\forall i  \label{1}
\end{equation}
The case $n_{6}\neq 0$ $\left( \text{and }d_{i}=0\text{ by }\left( 1\right)
\right) $ corresponds to those Lie algebras having non-abelian derived
subalgebra of dimension four.

\begin{enumerate}
\item  Suppose $n_{6}\neq 0$ and $e_{1}\neq 0$ $\left( \text{so }e_{j}=0%
\text{ for }i\geq 2\text{ }\right) .$\newline
We observe that if there exists an $a_{ij}\neq 0$ then the l.c.b. defined by 
$X_{2}^{\prime }=X_{2}+\alpha Y_{i}+\beta Y_{j}$ allows to suppose $%
f_{i}=f_{j}=0.$ So we have the conditions 
\begin{equation}
f_{i}a_{ij}=f_{j}a_{ij}=0,\;1\leq i,j  \label{2}
\end{equation}

\begin{enumerate}
\item  $\exists \,f_{i}\neq 0$ with $i\geq 2$ .We can suppose $f_{2}=1$ $%
\left( \text{so a}_{2j}=0\text{ by }\left( 2\right) \right) $ and $%
f_{i}=0,\forall i\geq 2.$

\begin{enumerate}
\item  If $a_{1j}=0$ for all $j$ we obtain two even dimensional algebras
isomorphic to $\frak{g}_{2m}^{24}$ and $\frak{g}_{2m}^{25}.$

\item  If $a_{1j}\neq 0$ for an index $j$ we can suppose $a_{13}=1.$ We
obtain two odd dimensional algebras isomorphic respectively to $\frak{g}%
_{2m+1}^{79}$ and $\frak{g}_{2m+1}^{80}.$
\end{enumerate}

\item  $f_{i}$ $=0,\;\forall i\geq 2$

\begin{enumerate}
\item  If $f_{1}\neq 0,$ then $a_{1j}=0$ by $\left( 2\right) .$ Reordering
the $\{Y_{2},..,Y_{n-6}\}$ we obtain an algebra isomorphic to $\frak{g}%
_{2m+1}^{81}.$

\item  If $f_{1}=0$ and $a_{12}\neq 0$ we obtain two algebras isomorphic to $%
\frak{g}_{2m}^{26}$ and $\frak{g}_{2m}^{27}.$

\item  If $f_{1}=a_{1j}=0,\;\forall j$ we obtain two algebras in odd
dimension isomorphic to $\frak{g}_{2m+1}^{82}$ and $\frak{g}_{2m+1}^{83}$.
\end{enumerate}
\end{enumerate}

\item  Suppose $n_{6}\neq 0$ and $e_{i}=0,\;\forall i.$\newline
We can suppose $p_{6}=0$ through a linear change of basis.

\begin{enumerate}
\item  If $f_{i}\neq 0$ for an index $i\geq 1$ let $f_{1}=1$ and $%
f_{i}=0,\;\forall i\geq 2$ and $a_{1j}=0$ by $\left( 2\right) .$ We obtain
two algebras isomorphic $\frak{g}_{2m+1}^{84}$ and $\frak{g}_{2m+1}^{85}.$

\item  If $f_{i}=0$ $\forall i$ we obtain two even dimensional algebras
respectively isomorphic to $\frak{g}_{2m}^{28}$ and $\frak{g}_{2m}^{29}.$
\end{enumerate}

\item  Suppose $n_{6}=0.$\newline
There are two cases to be separated: either there is a nonzero $d_{i}$ or
they are all zero. For simplicity, in the following we enumerate the cases
by fixing its parameters instead of deducing step by step as done up to
here. This method is first justified by the higher complexity of the
calculations in this case, for there are no conditions coming from the
characteristic sequence, and second because of the existence of linear
changes that would allow to jump from one step to another. These cases
complete the proof.

\begin{enumerate}
\item  $d_{1}\neq 0.$\newline
We can suppose $d_{1}=1$ and $d_{i}=0,\;i\geq 2.$ A change of type I allows
to suppose $o_{6}=0.$

\begin{enumerate}
\item  $d_{1}=e_{2}=f_{3}=1.$\newline
We have $e_{i}=0$ for $i\neq 2$ and $f_{i}=0$ for $i\neq 3.$ A change of
type I allows to suppose $p_{6}=0.$ Moreover, we have $a_{3j}=0,\;\forall
j\geq 4.$ Consider the change of Jordan basis defined by $X_{2}^{\prime
}=X_{2}+\alpha _{4}Y_{4}+\alpha _{5}Y_{5}+\alpha _{6}Y_{6}.$ If $a_{3j}=0$
for any index $i\geq 4,$ this change would allow to delete $f_{3},$ which is
not the case being considered. The remaining parameters can be structured as
follows
\bigskip
$
\begin{tabular}{l|l|l}
$\lbrack Y_{1},Y_{2}]=a_{12}X_{6}$ & $[Y_{1},Y_{4}]=\alpha X_{6},\;$ & $%
[Y_{2},Y_{4}]=a_{24}X_{6}$ \\ 
$\lbrack Y_{1},Y_{3}]=a_{13}X_{6}$ & and & $[Y_{2},Y_{5}]=\beta X_{6},$ \\ 
$\lbrack Y_{2},Y_{3}]=a_{23}X_{6}$ & $[Y_{1},Y_{i}]=0,\,i\geq 5$ & and \\ 
&  & $[Y_{2},Y_{i}]=0,\;i\geq 5$%
\end{tabular}
$ \newline
\newline
where $\alpha ,\beta =0,1$. If $a_{4i}\neq 0$ for $i>4$ or $a_{2j}\neq 0$
for $i=4$ or $i>5$ then we have $\alpha =0$ and $\beta =0$ respectively.%
\newline
We consider the distinct values for the pair $\left( a_{14},a_{25}\right)
=\left( \alpha ,\beta \right) :$\newline
A-1) $\left( \alpha ,\beta \right) =\left( 1,1\right) $\newline
A change of type III implies $a_{12}=a_{13}=a_{23}=a_{24}=0,$ so there
remains only an odd dimensional algebra isomorphic to $\frak{g}_{2m+1}^{86}.$%
\newline
A-2) $\left( \alpha ,\beta \right) =\left( 1,0\right) $\newline
A change of type III implies $a_{12}=a_{13}=0.$ The remaining parameters $%
a_{23}$ and $a_{24}$ are considered as a pair $\left( a_{23},a_{24}\right) :$
For $\left( 0,0\right) $ we obtain an algebra isomorphic to $\frak{g}%
_{2m}^{30},$ for $\left( 0,1\right) $ an algebra isomorphic to $\frak{g}%
_{2m}^{31},$ for $\left( 1,0\right) $ an algebra isomorphic to $\frak{g}%
_{2m}^{32}$ and for $\left( 1,1\right) $ one isomorphic to $\frak{g}%
_{2m}^{33}.$\newline
A-3) $\left( \alpha ,\beta \right) =\left( 0,1\right) $\newline
Type III implies $a_{12}=a_{23}=a_{24}=0.$ We obtain two algebras in even
dimension, which are respectively isomorphic to $\frak{g}_{2m}^{34}$ and $%
\frak{g}_{2m}^{35}.$\newline
A-4) $\left( \alpha ,\beta \right) =\left( 0,0\right) $\newline
If $a_{23}\neq 0$ we can suppose $a_{12}=0$ combining a change of type II
and a change of type IV. Moreover, combining type III and type IV we can
suppose $a_{13}=0,$ so we obtain a unique algebra in odd dimension and
isomorphic to $\frak{g}_{2m+1}^{87}.$\newline
If $a_{23}=0$ there are three algebras in odd dimension, which are
respectively isomorphic to $\frak{g}_{2m+1}^{88}$ if $\left(
a_{12},a_{13}\right) =\left( 0,0\right) ,\;\frak{g}_{2m+1}^{89}$ if $\left(
a_{12},a_{13}\right) =\left( 1,0\right) $ and $\frak{g}_{2m+1}^{90}\,$if $%
\left( a_{12},a_{13}\right) =\left( 0,1\right) .$ A combination of changes
type IV and II reduces the case $\left( a_{12},a_{13}\right) =(1,1)$ to the
last.

\item  $d_{1}=e_{2}=1,\;f_{i}=0$ for $i\geq 3.$\newline
We have $e_{i}=0$ for $i\neq 2.$ The remaining parameters are \newline
\newline
\begin{tabular}{l|l}
$\lbrack Y_{1},X_{2}]=X_{4}+f_{1}X_{6}$ & $[Y_{2},X_{2}]=X_{5}+f_{2}X_{6}$
\\ 
$\lbrack Y_{1},Y_{2}]=a_{12}X_{6}$ & $[Y_{2},Y_{4}]=\beta X_{6},\;\beta =0,1$
\\ 
$\lbrack Y_{1},Y_{3}]=\alpha X_{6},\;\alpha =0,1$ & and \\ 
$\lbrack Y_{2},Y_{3}]=a_{23}X_{6}$ & $[Y_{2},Y_{i}]=0$%
\end{tabular}
\newline
\newline
If $\alpha \neq 0$ then $a_{3j}=0,\;\forall j>3$ and if $\beta \neq 0$ then $%
a_{4i}=0$ for $i>4.$ Again we consider the pair $\left( \alpha ,\beta
\right) :$\newline
B-1) $\left( \alpha ,\beta \right) =\left( 1,1\right) $\newline
The change of Jordan basis defined by $X_{2}^{\prime }=X_{2}+AY_{3}+BY_{4}$
allows to suppose $f_{1}=f_{2}=0.$ A combination of changes type III and IV
implies $a_{12}=a_{23}=0,$ so we obtain a unique algebra isomorphic to $%
\frak{g}_{2m}^{36}.$\newline
B-2) $\left( \alpha ,\beta \right) =\left( 1,0\right) $\newline
A change of type IV allows to suppose $a_{12}=f_{1}=0.$ A change of type II
deletes $f_{2}$ and gives $\stackrel{\approx }{a_{13}}=1+a_{23}.$ Now, if $%
a_{23}\neq 0$ we consider the change $\frac{-1}{a_{23}}X_{1},..,\frac{1}{%
a_{23}^{2}}Y_{1},\frac{-1}{a_{23}^{3}}Y_{2}$ to obtain $a_{23}=-1,$ which
would imply $a_{13}=0$ by the previous change, but the assumption is $%
a_{13}=1,$ so it must be $a_{23}=0.$ We obtain a unique algebra in odd
dimension isomorphic to $\frak{g}_{2m+1}^{91}.$\newline
B-3) $\left( \alpha ,\beta \right) =\left( 0,1\right) $\newline
We have $a_{23}=0$ for the vector $Y_{3}$ is not distingushed any more from $%
Y_{j}$ for $j>3.$ If $a_{12}\neq 0$ a change of type II deletes $f_{1}$ and $%
f_{2}$ simultaneously. We obtain an algebra isomorphic to $\frak{g}%
_{2m}^{37}.$ If $a_{12}=0$ a similar change gives an algebra isomorphic to $%
\frak{g}_{2m}^{42}.$

\item  $d_{1}=1,\;e_{2}=0,\;f_{2}=1\;[o_{6}=f_{1}=0]$\newline
The remaining parameters are\newline
\newline
\begin{tabular}{l|l|l}
$\lbrack X_{3},X_{2}]=p_{6}X_{6}$ & $[Y_{1},Y_{3}]=\alpha X_{6},\;\alpha
=0,1 $ & $[Y_{2},Y_{3}]=a_{23}X_{6}$ \\ 
$\lbrack Y_{1},Y_{2}]=a_{12}X_{6}$ & and & $[Y_{2},Y_{4}]=0$ \\ 
& $[Y_{1},Y_{i}]=0$ for $i>4$ & 
\end{tabular}
\newline
\newline
where $a_{3j}=0$ for $j\geq 4$ if $\alpha \neq 0.$\newline
The parameter $a_{13}$ is zero, because the change $X_{2}^{\prime
}=X_{2}-Y_{3}$ would imply $f_{2}=0.$ In a similar way it must be $a_{23}=0.$
Depending on the values of $a_{12}$ and $\alpha $ we obtain four even
dimensional algebras which are respectively isomorphic to $\frak{g}%
_{2m}^{38},\frak{g}_{2m}^{39},\frak{g}_{2m}^{43}$ and $\frak{g}_{2m}^{44}.$

\item  $d_{1}=1,\;e_{i}=f_{i}=0,\;\forall i\geq 2.$\newline
With a change of type II we can suppose $e_{1}=0.$ \newline
D-1) If $a_{12}\neq 0$ a change of type I allows to write $f_{1}=0.$ We
obtain two algebras isomorphic to $\frak{g}_{2m}^{40}$ and $\frak{g}%
_{2m}^{41}.$\newline
D-2) If $a_{12}=0,$ then $f_{1}=0$ if $p_{6}\neq 0$ and we obtain an algebra
isomorphic to $\frak{g}_{2m+1}^{95}.$ If $p_{6}=0$ we obtain two algebras
respectively isomorphic to $\frak{g}_{2m+1}^{93}$ and $\frak{g}_{2m+1}^{94}.$
\end{enumerate}

\item  $d_{i}=0\;\forall i.$

\begin{enumerate}
\item  $e_{1}=f_{2}=1$\newline
We have $e_{i}=0$ for $i\geq 2$ and $f_{i}=0$ for $i\neq 2.$\newline
E-1) If there exists $j\geq 3$ with $a_{1j}\neq 0$ we can suppose $%
a_{13}=1,a_{1j}=0$ for $j>3.$ Moreover $a_{3j}=0$ for $j>3,$ for otherwise
it would be possible to delete $f_{2}$.\ Similarly we deduce $a_{2j}=0$ for $%
j>3.$ We obtain two odd dimensional algebras isomorphic to $\frak{g}%
_{2m+1}^{96}$ and $\frak{g}_{2m+1}^{97}.$\newline
E-2) If $a_{1j}=0$ for all $j$ then $a_{2j}=0$ like before and four
possibilities are given: in dependence on the values of the pair $\left(
a_{12},o_{6}\right) $ we obtain even dimensional algebras isomorphic
respectively to $\frak{g}_{2m}^{45},\frak{g}_{2m}^{46},\frak{g}_{2m}^{47}$
and $\frak{g}_{2m}^{48}.$

\item  $e_{1}=1,\;f_{i}=0$ $\forall i\geq 2.$\newline
F-1) If $a_{1j}\neq 0$ for $j\geq 2$ we suppose $a_{12}=1,\,a_{1j}=0$ for $%
j\neq 2.$ Then $a_{2j}=0$ for $j>2.$ A change of type I allows us to suppose 
$f_{1}=0.$ We obtain two algebras of even dimension and isomorphic to $\frak{%
g}_{2m}^{49}$ and $\frak{g}_{2m}^{50}.$ \newline
F-2) $a_{1j}=0$ for $j\geq 2.$ If $o_{6}\neq 0$ then a change of type II
allows to suppose $f_{1}=0.$ We obtain a unique algebra in odd dimension and
isomorphic to $\frak{g}_{2m+1}^{98}.$ For $o_{6}=0$ there are two algebras,
which are respectively isomorphic to $\frak{g}_{2m+1}^{99}$ and $\frak{g}%
_{2m+1}^{100}.$

\item  $e_{i}=0$ for $i\geq 1,$ $f_{1}\neq 0.$\newline
From this parameter values it follows $f_{i}=0$ for $i\geq 2$ and $a_{1j}=0$
for each $j$. We obtain three algebras respectively isomorphic to $\frak{g}%
_{2m+1}^{101},\;\frak{g}_{2m+1}^{102}$ and $\frak{g}_{2m+1}^{103}.$

\item  $e_{i}=f_{i}=0,\;\forall i.$\newline
A change of type II allows to suppose $p_{6}=0$ if $o_{6}\neq 0.$ We obtain
three even dimensional algebras isomorphic to $\frak{g}_{2m}^{51},\;\frak{g}%
_{2m}^{52}$ and $\frak{g}_{2m}^{53}.\;\;\;\;\;\;\;\;\;\;\;\;\;\;\;\;\;\;\;\;%
\;\;\;\;\;\;\;\;\;\;\;\;\;\;\;\;\;\;\;\;\;\;\;\;\;\;\;\;\;\;\;\;\;\;\;\;$%
\endproof%
\end{enumerate}
\end{enumerate}
\end{enumerate}
\end{enumerate}

\section{\textbf{Applications}}

We are now interested in those obtained laws which are characteristically
nilpotent, i.e, those of rank null. Characteristically nilpotent Lie
algebras were first introduced by Dixmier and Lister [9], and they have
become an important class of nilpotent algebras since then. Existence of
such algebras has been proved for any dimension $n\geq 7,$ as well as they
do not exist for $n\leq 6.$ There are a lot of papers constructing families
of characteristically nilpotent Lie algebras ( e.g [13], [14], [21]). As
most of known families and algebras are filiform, it is interesting to
obtain examples of characteristically nilpotent Lie algebras which are not
filiform. An interesting approach to this fact can be found in [14], where
characteristically nilpotent Lie algebras are obtained from the nilradical
of Borel subalgebras of complex simple Lie algebras. We will consider here
the $p$-filiform, characteristically nilpotent Lie algebras.

\begin{definition}
A Lie algebra $\frak{g}_{n}$ is called characteristically nilpotent if each
derivation $f\in Der\left( \frak{g}_{n}\right) $ is nilpotent.
\end{definition}

This definition results from a generalization of the descending central
sequence given by Dixmier and Lister. Unfortunately very little is known
about the algebra of derivations of a nilpotent Lie algebra, so that a
direct construction of a nilpotent Lie algebra of derivations is a not
trivial problem ( [9], [21].). However, characteristically nilpotent Lie
algebras behave as desired with sums, i.e, an algebra that is a finite sum
of ideals is characteristically nilpotent if and only if each ideal is
characteristically nilpotent [15]. From this we see that $\left( n-p\right) $%
-filiform characteristically nilpotent Lie algebras must be searched among
the nonsplit ones.

\begin{proposition}
For $p\leq 4$ there do not exist $\left( n-p\right) $-filiform
characteristically nilpotent Lie algebras.
\end{proposition}

\proof%
For the abelian and the Heisenberg algebras the assertion is evident. For $%
p=3$ and $4$ the proposition follows from the fact that these algebras have
all rank greater or equal than one ( see [6], [7]).%
\endproof%

\begin{remark}
It follows that characteristically nilpotent Lie algebras whose nilindex is
four must have characteristic sequence $\geq \left( 4,2,...,1\right) $.
\end{remark}

\begin{proposition}
A $\left( n-5\right) $-filiform Lie algebra is characteristically nilpotent
if and only if it is isomorphic to one of the following laws: $\frak{g}%
_{7}^{65},\frak{g}_{7}^{66,\alpha }$\newline
$\left( \alpha \neq 0\right) ,\frak{g}_{7}^{68},\frak{g}_{7}^{70},\frak{g}%
_{7}^{81},\frak{g}_{7}^{83},\frak{g}_{8}^{6},\frak{g}_{8}^{7,\alpha }\left(
\alpha \neq 0\right) ,\frak{g}_{8}^{9},\frak{g}_{8}^{11},\frak{g}_{8}^{25},%
\frak{g}_{8}^{27},\frak{g}_{9}^{57},\frak{g}_{9}^{80}$%
\endproof%
\end{proposition}

\begin{remark}
We observe that the seven dimensional family $\frak{g}_{7}^{66,\alpha }$ is
rigid in the variety $\frak{N}^{7}$ in the following sense: any perturbation
of a law in $\frak{g}_{7}^{66,\alpha }$ gives another law in this family.
Thus the closure of the orbit $\mathcal{O}\left( \frak{g}_{7}^{66,\alpha
}\right) $ in $\frak{N}^{7}$ gives an irreducible component of $\frak{N}^{7}$
[5].
\end{remark}

\begin{corollary}
There are characteristic nilpotent Lie algebras $\frak{g}_{n}$ with
nilpotence index $5$ for $n=7,8,9,14,15,16,17,18$ and $n\geq 21.$%
\endproof%
\end{corollary}

We have seen that the Lie algebra of derivations of a nilpotent Lie algebras
has not to be nilpotent in general. In fact the possibilities for the
algebra of derivations of a nilpotent Lie algebra are very ample. They can
vary from representations of the special linear algebra $\frak{sl}_{n}$ to
nilpotent Lie algebras, and no guide has been recognized until now. So it is
natural to ask for the existence of characteristically nilpotent Lie
algebras whose algebra of derivations has concrete properties: specifically
we ask if there are characteristically nilpotent Lie algebras of
derivations. That this doesn't always occur is shown by the following
example:

\begin{example}
Let $\frak{g}_{8}^{6}$ be the Lie algebra whose law is $\mu _{8}^{6}.$ The
algebra of derivations has dimension $13$ and is isomorphic to 
\[
\begin{tabular}{lll}
$\lbrack Z_{1},Z_{2}]=Z_{3},$ & $[Z_{2},Z_{3}]=-Z_{6},$ & $%
[Z_{3},Z_{10}]=-Z_{5}$ \\ 
$\left[ Z_{1},Z_{3}\right] =Z_{4},$ & $[Z_{2},Z_{6}]=-Z_{5},$ & $%
[Z_{3},Z_{13}]=-Z_{5}$ \\ 
$\lbrack Z_{1},Z_{4}]=Z_{5},$ & $[Z_{2},Z_{9}]=-Z_{6},$ & $%
[Z_{8},Z_{11}]=-Z_{5}$ \\ 
$\lbrack Z_{1},Z_{10}]=-Z_{6},$ & $[Z_{2},Z_{10}]=-Z_{6},$ & $%
[Z_{8},Z_{12}]=Z_{7}$ \\ 
$\lbrack Z_{1},Z_{11}]=-Z_{7},$ & $[Z_{2},Z_{12}]=Z_{5},$ & $%
[Z_{9},Z_{10}]=Z_{5}$ \\ 
& $[Z_{2},Z_{13}]=-Z_{4}$ & 
\end{tabular}
\;
\]
The linear system $\left( S\right) $ associated to this algebra has the
nontrivial solution 
\[
v=\left( \lambda _{i}\right) _{1\leq i\leq 13}=\lambda \left(
1,1,2,3,4,3,4,1,2,2,3,3,2\right) 
\]
so this algebra has nontrivial rank.
\end{example}

It seems that almost all caracteristically nilpotent Lie algebras will have
a non characteristically nilpotent Lie algebra of derivations.The existence
of algebras with characteristically nilpotent algebra of derivations is
proven by the next example, which gives a positive answer to the question
formulated by T\^{o}g\^{o} in [18]:

\begin{example}
For the algebra $\frak{g}_{7}^{81}$ with associated law $\mu _{7}^{81}$ the
algebra of derivations $Der\left( \frak{g}_{7}^{81}\right) $ has dimension $%
10$ and is isomorphic to the following algebra: 
\[
\text{%
\begin{tabular}{lll}
$\lbrack Z_{1},Z_{2}]=Z_{3},$ & $[Z_{2},Z_{6}]=-Z_{5},$ & $\left[
Z_{7},Z_{8}\right] =2Z_{5}-2Z_{6}+2Z_{10}$ \\ 
$\left[ Z_{1},Z_{3}\right] =Z_{4},$ & $[Z_{2},Z_{8}]=-Z_{6},$ & $%
[Z_{7},Z_{9}]=Z_{5}-2Z_{6}+2Z_{10}$ \\ 
$\lbrack Z_{1},Z_{4}]=Z_{5},$ & $[Z_{2},Z_{9}]=-Z_{4}-2Z_{6},$ & $%
[Z_{8},Z_{9}]=2Z_{6}-2Z_{10}$ \\ 
$\lbrack Z_{1},Z_{7}]=-Z_{4},$ & $[Z_{2},Z_{10}]=-Z_{5},$ &  \\ 
$\lbrack Z_{1},Z_{8}]=-Z_{6},$ & $[Z_{3},Z_{8}]=-Z_{5},$ &  \\ 
& $[Z_{3},Z_{9}]=-Z_{5},$ & 
\end{tabular}
} 
\]
It is not difficult to prove that this algebra is characteristically
nilpotent.
\end{example}

\begin{remark}
Thus it is possible to define an ''index'' for characteristically nilpotent
Lie algebras. The index equal to $1$ corresponds to the characteristically
nilpotent algebras like $\frak{g}_{8}^{6},$ i.e, those whose algebra of
derivations admits a nontrivial diagonalizable derivation. So we can call a
Lie algebra $\frak{g}$ characteristically nilpotent of index $k$ if $\frak{g}
$ and the $\left( k-1\right) $ first algebras of derivations are
characteristically nilpotent and the $k^{th}$ algebra of derivations is not
characteristically nilpotent. It would be interesting to know if there is a
relation between the nilpotence index or the characteristic sequence of the
algebra and the index $k$ defined above. It would be also interesting to
know if the sequence of derivation algebras stabilizes or not.
\end{remark}

\textit{\ }

\end{document}